\documentclass[11pt]{article}

\usepackage{color}
\usepackage{latexsym,dsfont}
\usepackage{amssymb}
\usepackage{amsmath}
\usepackage{amsthm}
\usepackage{graphicx}
\usepackage{enumerate}
\usepackage{hyperref}

\newtheorem{Theorem}{Theorem}[part]
\newtheorem{Definition}{Definition}[part]
\newtheorem{Proposition}{Proposition}[part]

\newtheorem{Lemma}{Lemma}[part]

\newtheorem{Corollary}{Corollary}[part]
\newtheorem{Remark}{Remark}[part]

\def \N{\mathbb{N}}
\def \R{\mathbb{R}}

\def \E{\mathbb{E}}
\def \F{\mathbb{F}}
\def \G{\mathbb{G}}
\def \H{\mathbb{H}}

\def \P{\mathbb{P}}
\def \Q{\mathbb{Q}}
\def \D{\mathbb{D}}

\def \Ac{{\cal A}}
\def \Bc{{\cal B}}
\def \Cc{{\cal C}}
\def \Dc{{\cal D}}
\def \Ec{{\cal E}}
\def \Fc{{\cal F}}
\def \Gc{{\cal G}}
\def \Hc{{\cal H}}

\def \Pc{{\cal P}}
\def \PMc{{\cal P \cal M}}
\def \Mc{{\cal M}}

\def \Sc{{\cal S}}

\def \Xc{{\cal X}}

\def \ni{\noindent}

\def \eps{\varepsilon}

\def \ep{\hbox{ }\hfill$\Box$}

\def\reff#1{{\rm(\ref{#1})}}

\def\beqs{\begin{eqnarray*}}
\def\enqs{\end{eqnarray*}}
\def\beq{\begin{eqnarray}}
\def\enq{\end{eqnarray}}

\newcommand{\nc}{\newcommand}
\nc{\esssup}{\mathop{\mathrm{ess\,sup}}}

\begin{document}

\title{Progressive enlargement of filtrations\\ and Backward SDEs with jumps\footnote{\textbf{Acknowledgement:} The authors would like to thank Shiqi Song for useful remarks which helped to improve the article.}
}
\author{Idris Kharroubi\footnote{The research of the author benefited from the support of the French ANR research grant LIQUIRISK.}\\ \footnotesize{CEREMADE, CNRS UMR 7534}, \\\footnotesize{Universit\'e Paris Dauphine}\\ 
\footnotesize{ \texttt{kharroubi @ ceremade.dauphine.fr}} 
  \and Thomas Lim\footnote{The research of the author benefited from the support of the ``Chaire Risque de Cr\' edit'', F\' ed\' eration Bancaire Fran\c caise.}\\ \footnotesize{Laboratoire d'Analyse et Probabilit\'es,}\\ \footnotesize{Universit\'e d'Evry}\footnotesize{ and ENSIIE,}\\ \footnotesize{\texttt{lim @ ensiie.fr}}
}

\date{June 2012}

\maketitle
\begin{abstract}
This work deals 
 with backward stochastic differential equation (BSDE) with random marked jumps, and their applications to default risk. We show that these BSDEs are linked with Brownian BSDEs through the decomposition of processes 
  with respect to the progressive enlargement of filtrations. We prove that the equations have solutions if the associated Brownian BSDEs have solutions. We also provide a uniqueness theorem for BSDEs with jumps by giving a comparison theorem based on the comparison for Brownian BSDEs. We give in particular some results  for quadratic BSDEs. As applications, we study the pricing and the hedging of a European option in a market with a single jump, and the utility maximization problem in an incomplete market with a finite number of jumps.
\end{abstract}

\vspace{1cm}


\noindent \textbf{MSC classification (2000):} 60G57, 60J75, 91G10, 93E20.


\setcounter{section}{0}


\section{Introduction}
\setcounter{equation}{0} \setcounter{Assumption}{0}
\setcounter{Theorem}{0} \setcounter{Proposition}{0}
\setcounter{Corollary}{0} \setcounter{Lemma}{0}
\setcounter{Definition}{0} \setcounter{Remark}{0}

In recent years, credit risk has come out to be one of most fundamental financial risk. The most extensively studied form of credit risk is the default risk. Many people, such as Bielecki, Jarrow, Jeanblanc, Pham, Rutkowski (\cite{bierut04, biejearut04a, jaryu01, jealec09, jiapha09, pha09}) and many others, have worked on this subject. In several papers (see for example Ankirchner \emph{et al.} \cite{ankblaeyr09}, Bielecki and Jeanblanc \cite{biejea08} and  Lim and Quenez \cite{limque09}), related to this topic, backward stochastic differential equations (BSDEs) with jumps have appeared. Unfortunately, the results relative to these latter BSDEs are far from being as numerous as for Brownian BSDEs. In particular, there is not any general result on the existence and the uniqueness of solution to quadratic BSDEs, except Ankirchner \emph{et al.} \cite{ankblaeyr09}, in which the assumptions on the driver are strong. In this paper, we study BSDEs with random marked jumps and apply the obtained results to mathematical finance where these jumps can be interpreted as default times. We give a general existence and uniqueness result for the solutions to these BSDEs, in particular we enlarge the result given by \cite{ankblaeyr09} for quadratic BSDEs. 

A standard approach of credit risk modeling is based on the powerful technique of filtration enlargement, by making the distinction between the filtration $\F$ generated by the Brownian motion, and its smallest extension $\G$ that turns default times into $\G$-stopping times. This kind of filtration enlargement has been referred to as progressive enlargement of filtrations. This field is a traditional subject in probability theory initiated by fundamental works of the French school in the 80s, see e.g. Jeulin \cite{jeu80}, Jeulin and Yor \cite{jeuyor85}, and Jacod \cite{jac87}. For an overview of applications of progressive enlargement of filtrations on credit risk, we refer to the books of Duffie and Singleton \cite{dufsin03}, of Bielecki and Rutkowski \cite{bierut04}, or the lectures notes of Bielecki \emph{et al.} \cite{biejearut04a}. 

The purpose of this paper is to combine results on Brownian BSDEs and  results on progressive enlargement of filtrations in view of providing existence and uniqueness of solutions to BSDEs with random marked jumps. We consider a progressive enlargement with multiple random times and associated marks. These marks can represent for example the name of the firm which defaults or the jump sizes of asset values. 

Our approach consists in using the recent results of Pham \cite{pha09} on the decomposition of predictable processes with respect to the progressive enlargement of filtrations to decompose a BSDE with random marked jumps into a sequence of Brownian BSDEs. By combining the solutions of Brownian BSDEs, we obtain a solution to the BSDE with random marked times. This method allows to get a general existence theorem. In particular, we get an existence result for quadratic BSDEs which is more general than the result of Ankirchner \textit{et al} \cite{ankblaeyr09}. This decomposition approach also allows to obtain a uniqueness theorem under  Assumption \textbf{(H)} i.e. 
 any $\F$-martingale remains a $\G$-martingale. 
 We first  set a general comparison theorem  for BSDEs with jumps based on  comparison theorems for Brownian BSDEs. Using this theorem, we prove, in particular, the uniqueness for quadratic BSDEs with a concave generator w.r.t. $z$.

We illustrate our methodology with two financial applications in default risk management: the pricing and the hedging of a European option, and the problem of utility maximization in an incomplete market. A similar problem (without marks) has recently been considered in Ankirchner \emph{et al.} \cite{ankblaeyr09} and Lim and Quenez \cite{limque09}.

The paper is organized as follows. The next section presents the general framework of progressive enlargement of filtrations with successive random times and marks, and states the decomposition result for $\G$-predictable and  specific $\G$-progressively measurable processes. In Section 3, we use this decomposition to make a link between Brownian BSDEs and BSDEs with random marked jumps. This allows to give a general existence result under a density assumption. We then give two examples: quadratic BSDEs with marked jumps for the first one, and linear BSDEs arising in the pricing and hedging problem of a European option in a market with a single jump for the second one. In Section 4, we give a general comparison theorem for BSDEs and we use this result to give a uniqueness theorem for quadratic BSDEs. Finally, in Section 5, we apply our  existence and uniqueness results  to solve the exponential utility maximization problem in an incomplete market with a finite number of marked  jumps.

\section[Progressive enlargement of filtrations]{Progressive enlargement of filtrations with successive random times and marks}
\setcounter{equation}{0} \setcounter{Assumption}{0}
\setcounter{Theorem}{0} \setcounter{Proposition}{0}
\setcounter{Corollary}{0} \setcounter{Lemma}{0}
\setcounter{Definition}{0} \setcounter{Remark}{0}

We fix a probability space $(\Omega, \Gc, \P)$, and we start with a reference filtration $\F=(\Fc_{t})_{t \geq 0}$ satisfying the usual conditions\footnote{$\Fc_{0}$ contains the $\P$-null sets and $\F$ is right continuous: $\Fc_{t}=\Fc_{t^+}:=\cap_{s>t}\Fc_{s}$.} and generated by a $d$-dimensional Brownian motion $W$. Throughout the sequel, we consider a finite sequence $(\tau_{k}, \zeta_{k})_{1 \leq k\leq n}$, where 
\begin{itemize}
\item $(\tau_{k})_{1 \leq k\leq n}$ is a nondecreasing sequence of random times (i.e. nonnegative $\Gc$-random variables), 
\item  $(\zeta_{k})_{1 \leq k\leq n}$ is a sequence of random marks valued in some Borel subset $E$ of $\R^m$. 
\end{itemize}
We denote by $\mu$ the random measure associated with the sequence $(\tau_{k}, \zeta_{k})_{1\leq k\leq n}$ : 
\beqs
\mu([0,t] \times B) & = &\sum_{k=1}^n\mathds{1}_{\{\tau_{k} \leq t,\,\zeta_{k}\in B\}}\;,\quad t\geq 0\;, ~B\in\Bc(E)  \;.
\enqs

For each $k=1,\ldots,n$, we consider $\D^k=(\Dc^k_{t})_{t \geq 0}$ the smallest filtration for which $\tau_{k}$ is a stopping time and $\zeta_{k}$ is $\Dc_{\tau_{k}}^k$-measurable. $\D^k$ is then given by 
$ \Dc _{t}^k=\sigma(\mathds{1}_{\tau_{k}\leq s}, s\leq t)\vee\sigma(\zeta_{k}\mathds{1}_{\tau_{k}\leq s}, s\leq t)$.  The global information is then defined by the progressive enlargement $\G=(\Gc_{t})_{t \geq 0}$ of the initial filtration $\F$ where $\G$ is the smallest right-continuous filtration containing $\F$, and such that for each $k=1,\ldots,n$, $\tau_{k}$ is a $\G$-stopping time, and $\zeta_{k}$ is $\Gc_{\tau_{k}}$-measurable. $\G$ is given by $\Gc_t=\tilde \Gc_{t^+}$, where $\tilde \Gc _t = \Fc_t\vee \Dc^1_t\vee \cdots \vee \Dc^n_t$ for all $t\geq 0$.

We denote by $\Delta_k$ the set where the random $k$-tuple $(\tau_1,\ldots,\tau_k)$ takes its values in $\{\tau_n<\infty\}$:
\beqs
\Delta_{k} & := & \big\{(\theta_{1},\ldots,\theta_{k})\in(\R_{+})^k~:~\theta_{1}\leq\ldots\leq \theta_{k}\big\}, \quad  1 \leq k\leq n\;.
\enqs
We introduce some notations used throughout the paper:
\begin{itemize}
\item $\Pc(\F)$ (resp. $\Pc(\G)$) is the $\sigma$-algebra of $\F$ (resp. $\G$)-predictable measurable subsets of $\Omega \times \R_{+}$, i.e. the $\sigma$-algebra generated by the left-continuous $\F$ (resp. $\G$)-adapted processes. 
\item $\PMc(\F)$ (resp. $\PMc(\G)$) is the $\sigma$-algebra of $\F$ (resp. $\G$)-progressively measurable subsets of $\Omega\times\R_{+}$. 
\item For $k=1,\ldots,n$, $\PMc(\F,\Delta_k,E^k)$ 
is the $\sigma$-algebra  generated by processes $X$ from $\R_+\times\Omega\times\Delta_k\times E^k$ to $\R$ such that  $(X_t(.))_{t\in[0,s]}$ is $\Fc_s\otimes\Bc([0,s]) \otimes\Bc(\Delta_{k})\otimes\Bc(E^k)$-measurable, for all $s\geq 0$. 
\item For $\theta = (\theta_{1},\ldots,\theta_n)\in\Delta_{n}$ and $e=(e_{1},\ldots,e_n) \in E^n$, we denote by
\beqs
\theta_{(k)} &= &(\theta_{1},\ldots, \theta_{k})~~\mbox{ and }~~ e_{(k)} ~= ~(e_{1},\ldots, e_{k})\;,\quad  1 \leq k\leq n\;.
\enqs
We also denote by $\tau_{(k)}$ for $(\tau_{1},\ldots, \tau_{k})$ and $\zeta_{(k)}$ for $(\zeta_{1},\ldots, \zeta_{k})$, for all $k = 1, \ldots, n$.
\end{itemize}

The following result provides the basic decomposition of predictable and progressive processes with respect to this progressive enlargement of filtrations.

\begin{Lemma}\label{3lemme decomposition}
\begin{enumerate}[(i)]

\item Any $\Pc(\G)$-measurable process $X=(X_t)_{t \geq 0}$ is represented as
\beq
\label{3decomposition previsible}
X_t & = & X^0_t\mathds{1}_{t\leq \tau_1}+\sum_{k=1}^{n-1}X^k_t( \tau_{(k)},\zeta_{(k)})\mathds{1}_{\tau_k<t\leq \tau_{k+1}}
+X^n_t(\tau_{(n)},\zeta_{(n)})\mathds{1}_{\tau_n<t}\;,
\enq
for all $t \geq 0$, where $X^0$ is $\mathcal{P}(\F)$-measurable and $X^k$ is $\mathcal{P}(\F) \otimes\Bc(\Delta_{k})\otimes\Bc(E^k)$-measurable for $k=1,\ldots,n$. 
\item  Any c\`{a}d-l\`{a}g $\PMc(\G)$-measurable process $X=(X_t)_{t \geq 0}$ of the form
\beqs
X_t & = & J_t+\int_0^t\int_EU_s(e)\mu(de,ds)\;, \quad t\geq 0\;,
\enqs
where $J$ is $\Pc(\G)$-measurable and $U$ is $\Pc(\G)\otimes\Bc(E)$-measurable,
is represented as
\beq
\label{3decomposition progressive}
X_t & = & X^0_t\mathds{1}_{t < \tau_1}+\sum_{k=1}^{n-1}X^k_t( \tau_{(k)},\zeta_{(k)})\mathds{1}_{\tau_k \leq t < \tau_{k+1}}
+X^n_t(\tau_{(n)},\zeta_{(n)})\mathds{1}_{\tau_n \leq t}\;,
\enq
for all $t \geq 0$, where $X^0$ is $\PMc(\F)$-measurable and $X^k$ is $\PMc(\F,\Delta_k,E^k)$-measurable for $k=1,\ldots,n$.
\end{enumerate}
\end{Lemma}

The proof of (i) is given in Pham \cite{pha09} and is therefore omitted. The proof of (ii) 
 is based on similar arguments. Hence, we postpone it to the appendix.\\

Throughout the sequel, we will use the convention $\tau_0=0$, $\tau_{n+1}=+\infty$, $\theta_0=0$ and $\theta_{n+1}=+\infty$ for any $\theta\in\Delta_n$, and $X^0(\theta_{(0)},e_{(0)})= X^0$ to simplify the notation. 
\\

\begin{Remark}\label{remparam}
%
 {\rm In the case where the studied process $X$ depends on another parameter $x$ evolving in a Borelian subset $\mathcal{X}$ of $\R^p$, and if $X$ is $\Pc(\G)\otimes \Bc(\mathcal{X})$, 
  then, decomposition \reff{3decomposition previsible} 
   is still true but where $X^k$ is $\Pc(\F)\otimes\Bc(\Delta_k)\otimes\Bc(E^k)\otimes\Bc(\mathcal{X})$
   -measurable. Indeed, it is obvious for the processes generating $\Pc(\G)\otimes\Bc(\mathcal{X})$ 
   of the form $X_t(\omega, x)= L_t(\omega)R(x)$, $(t,\omega,x)$ $\in$ $\R_+\times\Omega\times\mathcal{X}$, where $L$ is $\Pc(\G)$
   -measurable and $R$ is $\Bc(\mathcal{X})$-measurable. Then, the result is extended to any $\Pc(\G)\otimes\Bc(\mathcal{X})$
   -measurable process by the monotone class theorem.
 }
\end{Remark}
We now introduce a density assumption on the random times and their associated marks by assuming that  
the distribution of $(\tau_{1},\ldots,\tau_{n},\zeta_{1},\ldots,\zeta_{n})$  
is absolutely continuous with respect to the Lebesgue measure $d\theta \,de$ on $\Bc(\Delta_{n})\otimes\Bc(E^n)$. More precisely, we make the following assumption.

\vspace{3mm}

\ni (HD) There exists a positive  $\Pc(\F)\otimes\Bc(\Delta_{n})\otimes\Bc( E^n)$-measurable map $\gamma$  
such that for any $t\geq 0$, 
\beqs
\P[(\tau_{1},\ldots,\tau_{n}, \zeta_1, \ldots, \zeta_n)\in d\theta de |\Fc_t]  =  \gamma_t
(\theta_{1},\ldots,\theta_{n},e_{1},\ldots,e_{n}) d\theta_{1}\ldots d\theta_{n} de_1\ldots de_n\;.
\enqs

\vspace{3mm}



\ni We then introduce some notation. Define the process $\gamma^0$ by
\beqs
\gamma^0_t & = & \P[\tau_1>t|\Fc_t]
 ~ = ~ \int_{\Delta_n\times E^n}\mathds{1}_{\theta_1>t } \gamma_t(\theta,e)d\theta de\;,
 \enqs
 and the map $\gamma^k$ a $\Pc(\F)\otimes\Bc(\Delta_k)\otimes\Bc(E^k)$-measurable process, $k$ $=$ $1,\ldots,n-1$, by
 \beqs
& &  \gamma^{k}_t\big(\theta_{1},\ldots,\theta_{k},e_1,\ldots,e_{k}\big) \\
& = &  \int_{\Delta_{n-k}\times E^{n-k}}\mathds{1}_{\theta_{k+1}>t}\gamma_t(\theta_1,\ldots,\theta_n,e_1\ldots,e_n)d\theta_{k+1}\ldots d\theta_n de_{k+1}\ldots de_n\;.
\enqs
We shall use the natural convention $\gamma^n=\gamma$.
We obtain that under (HD), the random measure $\mu$ admits a compensator absolutely continuous w.r.t. the Lebesgue measure. The intensity $\lambda$ is given by the following proposition. 
\begin{Proposition} \label{prop intensite}
Under (HD), the random measure $\mu$ admits a compensator for the filtration $\G$ given by $\lambda_t(e)dedt$,
where the intensity $\lambda$ is defined by
\beq\label{decomp lambda}
\lambda_t(e) & = & 
\sum_{k=1}^n\lambda^k_t(e,\tau_{(k-1)},\zeta_{(k-1)})\mathds{1}_{\tau_{k-1}< t \leq \tau_{k}}\;,
\enq
with 
\beqs
\lambda_t^k(e,\theta_{(k-1)},e_{(k-1)}) & = & \frac{\gamma_t^{k}(\theta_{(k-1)},t,e_{(k-1)},e)}{\gamma^{k-1}_t(\theta_{(k-1)},e_{(k-1)})}\;, \quad(\theta_{(k-1)},t,e_{(k-1)},e)\in \Delta_{k-1} \times \R_+\times E^{k}\;.
\enqs
\end{Proposition}
The proof of Proposition \ref{prop intensite} is based on similar arguments to those of \cite{nekjeajia}. We therefore postpone it to the appendix.\\

We add an assumption on the intensity $\lambda$ which will be used in existence and uniqueness results for quadratic BSDEs as well as for the utility maximization problem:
\beqs
\text{(HBI)}\qquad\qquad\text{The process}\quad \Big(\int_E \lambda_t(e) de \Big)_{t \geq 0} \quad \text{is bounded on } [0, \infty) \;.
\enqs

We now consider one dimensional BSDEs driven by $W$ and the random measure $\mu$. To define solutions, we need to introduce the following spaces, where $a,b\in\R_+$ with $a\leq b$,
and $T < \infty$ is the terminal time:

\begin{itemize}
\item $\Sc_{\G}^\infty[a,b]$ (resp. $\Sc_{\F}^\infty[a,b]$) is the set of $\R$-valued $\PMc(\G)$ (resp.  $\PMc(\F)$)-measurable processes $(Y_t)_{t\in[a,b]}$ essentially bounded:
\beqs
{\| Y \|}_{\Sc^\infty[a,b]} & := & \esssup_{t\in[a,b]}|Y_{t}|~<~\infty \;.
\enqs

\item $L^2_{\G}[a,b]$ (resp. $L^2_{\F}[a,b]$) is the set of $\R^d$-valued $\Pc(\G)$ (resp. $\Pc(\F)$)-measurable processes $(Z_t)_{t\in[a,b]}$ such that
\beqs
\|Z\|_{L^2[a,b]}  & := & \Big(\E\Big[ \int_a^b |Z_t|^2 dt \Big]\Big)^{1\over2} ~< ~\infty \;.
\enqs
\item $L^2(\mu)$ is the set of  $\R$-valued $\Pc(\G)\otimes\Bc(E)$-measurable processes $U$ such that
\beqs
\|U\|_{L^2(\mu)} &:= & \Big(\E\Big[\int_0^T\int_{E}|U_{s}(e)|^2\mu(de,ds)\Big]\Big)^{1\over2}~<~\infty\;.
\enqs

\end{itemize}

\ni We then consider BSDEs of the form: find a triple  $(Y,Z,U)$ $\in$ $\Sc^\infty_\G[0,T] \times L^2_\G[0,T] \times L^2(\mu)$ such that\footnote{The symbol $\int_{s}^t$ stands for the integral on the interval $(s,t]$ for all $s,t$ $\in$ $\R_{+}$.}
\beq\label{3BSDE jump}
Y_t ~  = ~  \xi+\int_t^Tf(s,Y_s,Z_s,U_s)ds-\int_t^TZ_sdW_s-\int_t^T\int_{E}U_s(e)\mu(de,ds),~0\leq t\leq T,
\enq
where 
\begin{itemize}
\item $\xi$  is a $\Gc_T$-measurable random variable of the form: 
\beq \label{3decomposition xi}
\xi & = &\sum_{k=0}^{n}\xi^k(\tau_{(k)}, \zeta_{(k)})\mathds{1}_{\tau_k\leq T< \tau_{k+1}} \;,
\enq
with $\xi^0$ is $\Fc_{T}$-measurable and $\xi^k$ is $\Fc_{T}\otimes\Bc(\Delta_{k})\otimes\Bc(E^{k})$-measurable for each $k=1,\ldots,n$, 
\item  $f$ is map from $[0,T]\times\Omega\times\R\times\R^d\times Bor(E,\R)$ to $\R$ which is  a $\Pc(\G)\otimes\Bc(\R)\otimes\Bc(\R^d)\otimes\Bc(Bor(E,\R))$-$\Bc(\R)$-measurable map. 
Here,  $Bor(E,\R)$ is the set of borelian functions from $E$ to $\R$, and  $\Bc(Bor(E,\R))$ is the borelian $\sigma$-algebra on $Bor(E,\R)$ for the pointwise convergence topology. 
\end{itemize}

\ni To ensure that BSDE \reff{3BSDE jump} is well posed, we have to check that the stochastic integral w.r.t. $W$ is well defined on $L^2_\G[0,T]$ in our context. 
\begin{Proposition}
Under (HD), for any process $Z\in L^2_\G[0,T]$, the stochastic integral $\int_0^TZ_sdW_s$ is well defined.
\end{Proposition}
\ni \textbf{Proof.}
Consider the initial progressive enlargement $\H$ of the filtration $\G$. We recall that $\H=(\Hc_t)_{t\geq 0}$ is given by 
\beqs
\Hc_t & = & \Fc_t\vee\sigma\big( \tau_1,\ldots,\tau_n,\zeta_1,\ldots,\zeta_n \big)\;,\quad t\geq 0\;.
\enqs
We prove that the stochastic integral $\int_0^TZ_sdW_s$ is well defined for all $\Pc(\H)$-measurable process $Z$ such that $\E\int_0^T|Z_s|^2ds<\infty$.
Fix such a process $Z$.

From Theorem 2.1 in \cite{jac87}, we obtain that $W$ is an $\H$-semimartingale of the form
\beqs
W_t & = & M_t+\int_0^ta_s(\tau_{(n)},\zeta_{(n)})ds\;,\quad t\geq 0\;,
\enqs
where $a$ is $\Pc(\F)\otimes\Bc(\Delta_n)\otimes\Bc(E^n)$-measurable. Since $M$ is a $\H$-local continuous martingale with quadratic variation $\langle M, M \rangle_t=\langle W, W \rangle_t = t$ for $t\geq 0$, we get from L\'evy's characterization of Brownian motion (see e.g. Theorem 39 in \cite{pro05}) that $M$ is a $\H-$Brownian motion. 
Therefore the stochastic integral $\int_0^TZ_sdM_s$ is well defined  and we now concentrate on the term $\int_0^T Z_sa_s(\tau_{(n)},\zeta_{(n)})ds$.

From Lemma 1.8 in \cite{jac87} the process $\gamma(\theta,e)$ is an $\F$-martingale. Since $\F$ is the filtration generated by $W$ we get from the representation theorem  of Brownian martingales that 
\beqs
\gamma_t(\theta,e) & = & \gamma_0(\theta,e)+\int_0^t\Gamma_s(\theta,e)dW_s\;,\quad t\geq 0\;.
\enqs
Still using Theorem 2.1 in \cite{jac87} and since $\gamma(\theta,e)$ is continuous, we have 
\beqs
\langle \gamma(\theta,e),W\rangle_t & = & \int_0^t\gamma_{s}(\theta,e)a_s(\theta,e)ds\;,\quad t\geq 0 
\enqs
for all $(\theta,e)\in \Delta_n\times E^n$. Therefore we get 
\beqs
\Gamma_s(\theta,e) & = & \gamma_s(\theta,e) a_s(\theta,e)\;,\quad s\geq 1
\enqs
for all $(\theta,e)\in \Delta_n\times E^n$. Since $\gamma(\theta,e)$ is an $\F$-martingale, we obtain (see e.g.  Theorem 62 Chapter 8 in \cite{demey80}) that
\beq\label{eq merdique}
\int_0^T|\gamma_s(\theta,e)a_s(\theta,e)|^2ds & < & +\infty\;, \quad \P-a.s.
\enq
for all $(\theta,e)\in\Delta_n\times E^n$. 
Consider the set $A\in\Fc_T\otimes\Bc(\Delta_n)\otimes\Bc(E^n)$ defined by
\beqs
A & := & \Big\{(\omega,\theta,e)\in\Omega\times\Delta_n\times E^n~:~\int_0^T|\gamma_s(\theta,e)a_s(\theta,e)|^2ds = +\infty\;, 
\Big\}\;.
\enqs
Then, we have $\P(\tilde\Omega)=0$, where 
\beqs
\tilde \Omega & = & \big\{\omega\in\Omega~:~\big(\omega,\tau(\omega),\zeta(\omega)\big)\in A \big\} \;.
\enqs
Indeed, we have from the density assumption (HD)
\beq\nonumber
\P(\tilde\Omega) & = & \E\Big[\mathds{1}_{A}\big(\omega,\tau(\omega),\zeta(\omega)\big)\Big]~=~ \E\Big[\E\Big[\mathds{1}_{A}\big(\omega,\tau(\omega),\zeta(\omega)\big)\Big|\Fc_T\Big]\Big]\\
 & = & \int_{\Delta_n\times E^n}\E\Big[\mathds{1}_{A}\big(\omega,\theta,e\big)\gamma_T(\theta,e)\Big]d\theta de\;.\label{eq p omega tilde }
\enq
From the definition of $A$ and \reff{eq merdique}, we have 
\beqs
\mathds{1}_{A}\big(.,\theta,e\big)\gamma_T(\theta,e) & = & 0 \;,\quad \P-a.s.
\enqs
for all $(\theta,e)\in \Delta_n\times E^n$.
Therefore, we get from \reff{eq p omega tilde }, $\P(\tilde\Omega)=0$ or equivalently
\beq\label{int ga care fini}
\int_0^T \big|\gamma_s(\tau_1,\ldots,\tau_n,\zeta_1,\ldots,\zeta_n)a_s(\tau_1,\ldots,\tau_n,\zeta_1,\ldots,\zeta_n)\big|^2ds & < & +\infty\;, \quad \P-a.s.
\enq
From Corollary 1.11 we have $\gamma_t(\tau_1,\ldots,\tau_n,\zeta_1,\ldots,\zeta_n)>0$ for all $t\geq 0$ $\P$-a.s. Since $\gamma_.(\tau_1,\ldots,\tau_n,\zeta_1,\ldots,\zeta_n)$ is continuous we obtain
\beq\label{gamma ne s'annule pas}
\inf_{s\in[0,T]} \gamma_s(\tau_1,\ldots,\tau_n,\zeta_1,\ldots,\zeta_n) & > & 0\;, \quad \P-a.s.
\enq
Combining \reff{int ga care fini} and \reff{gamma ne s'annule pas}, we get 
\beqs
\int_0^T \big| a_s(\tau_1,\ldots,\tau_n,\zeta_1,\ldots,\zeta_n) \big|^2ds & < & +\infty\;, \quad \P-a.s.
\enqs
 Since $Z$ satisfies  $\E\int_0^T |Z_s|^2ds<\infty$, we obtain that 
 \beqs
 \int_0^T \big|Z_sa_s(\tau_1,\ldots,\tau_n,\zeta_1,\ldots,\zeta_n) \big|ds & < & +\infty\;,\quad \P-a.s.
 \enqs
 Therefore  $\int_0^T Z_sa_s(\tau_1,\ldots,\tau_n,\zeta_1,\ldots,\zeta_n)ds$ is well defined.  
\ep
\section{Existence of a solution}
\setcounter{equation}{0} \setcounter{Assumption}{0}
\setcounter{Theorem}{0} \setcounter{Proposition}{0}
\setcounter{Corollary}{0} \setcounter{Lemma}{0}
\setcounter{Definition}{0} \setcounter{Remark}{0}

In this section, we use the decompositions given by Lemma \ref{3lemme decomposition} to solve BSDEs with a finite number of jumps. We use a similar approach to Ankirchner \emph{et al.} \cite{ankblaeyr09}: one can explicitly construct a solution by combining solutions of an associated recursive system of Brownian BSDEs. But contrary to them, we suppose that there exist $n$ random times and $n$ random marks. Our assumptions on the driver are also weaker.  
Through a simple example we first show how our method to construct solutions to BSDEs with jumps works. We then give  a general existence theorem which links the studied BSDEs with jumps with a system of recursive Brownian BSDEs. We finally illustrate our general result with concrete examples. 

\subsection{An introductory example}
We begin by giving a simple example to illustrate the used method. We consider the following equation involving only a single jump time $\tau$ and a single mark $\zeta$ valued in $E=\{0,1\}$: 
\begin{equation}\label{EDSR exemple}
\left\{\begin{array}{rcl}Y_T &  = &  c \mathds{1}_{T < \tau} + h(\tau,\zeta) \mathds{1}_{T \geq \tau}\\
-dY_t & = & f(U_t) dt - U_t dH_t \;,\quad 0\leq t\leq T\;,
\end{array}\right.
\end{equation}
where $H_t =(H_t(0),H_t(1))$ with $H_t(i)=\mathds{1}_{\tau\leq t,\zeta=i}$ for $t\geq0$ and $i\in E$. Here $c$ is a real constant, and $f$ and $h$ are deterministic functions. To solve BSDE \reff{EDSR exemple}, we first solve a recursive system of BSDEs:
\beqs
Y^1_t (\theta,e) & = & h(\theta,e) + f(0,0) (T-t) \;, \quad \theta \wedge T \leq t \leq T \;, \\
Y^0_t & = & c + \int_t^T f \big( Y^1_s(s,0) - Y^0_s, Y^1_s(s,1) - Y^0_s \big) ds \;, \quad 0 \leq t \leq T \;.
\enqs 
Suppose that the recursive system of BSDEs admits for any $(\theta,e) \in [0,T] \times \{0,1\}$ a couple of solution $Y^1(\theta,e)$ and $Y^0$.
Define the process $(Y, U)$ by 
\beqs
Y_t &  =  & Y^0_t \mathds{1}_{t < \tau} + Y^1_t(\tau,\zeta) \mathds{1}_{t \geq \tau}\;,\quad t\in [0,T]\;,
\enqs
and
\beqs
U_t(i) = (Y^1_t(t,i) -Y^0_t) \mathds{1}_{t \leq \tau}\;,\quad t\in [0,T]\;,~i=0,1\;.
\enqs 
We then prove that the process $(Y, U)$ is solution of BSDE \reff{EDSR exemple}. By It\^o's formula, we have 
\beqs
dY_t & = & d\Big( Y^0_t \mathds{1}_{t < \tau} + Y^1_t(\tau,\zeta) \mathds{1}_{t \geq \tau} \Big)\\
 & =  & d\Big( Y^0_t (1-H_t(0)-H_t(1)) +\int_0^th(s,0)dH_s(0)\\
  & &+\int_0^th(s,1)dH_s(1) + (H_t(0)+H_t(1)) f(0,0)(T-t) \Big)\;. 
\enqs
This can be written
\beqs
dY_t & = & - \big[(1-H_t(0)-H_t(1)) f \big( Y^1_t(t,0) - Y^0_t, Y^1_t(t,1) - Y^0_t \big) + (H_t(0)+H_t(1)) f(0,0) \big] dt  \\
 & & + \big[ h(t,0) + (T-t) f(0,0) - Y^0_t \big] dH_t(0)+ \big[ h(t,1) + (T-t) f(0,0) - Y^0_t \big] dH_t(1) \;. 
\enqs
From the definition of $U$, we get 
\beqs
dY_t & = & - f(U_t) dt + U_t dH_t \;.  
\enqs
We also have $Y_T =  c \mathds{1}_{T < \tau} + h(\tau,\zeta) \mathds{1}_{T \geq \tau}$, which shows that  $(Y, U)$ is solution of BSDE \reff{EDSR exemple}. 
\subsection{The existence theorem}

To prove the existence of a solution to BSDE \reff{3BSDE jump}, we introduce the decomposition of the coefficients $\xi$ and $f$ as given by \reff{3decomposition xi} and Lemma \ref{3lemme decomposition}. 

From   Lemma \ref{3lemme decomposition} (i) and Remark \ref{remparam}, we get the following decomposition for $f$
\beq \label{3decomposition f}
f(t,y,z,u) & =  & \sum_{k=0}^{n}f^k(t, y, z, u, \tau_{(k)},  \zeta_{(k)})
\mathds{1}_{\tau_k\leq t< \tau_{k+1}} \;,
\enq
where $f^0$ is $\Pc(\F)\otimes\Bc(\R)\otimes\Bc(\R^d)\otimes\Bc(Bor(E,\R))$-measurable and $f^k$ is $\Pc(\F)\otimes\Bc(\R)\otimes\Bc(\R^d)\otimes\Bc(Bor(E,\R))\otimes\Bc(\Delta_{k})\otimes\Bc(E^{k})$-measurable for each $k=1,\ldots,n$.

In the following theorem, we show how BSDEs driven by $W$ and $\mu$ are related to a recursive system of Brownian BSDEs involving the coefficients $\xi^k$ and $f^k$, $k=0,\ldots,n$. 
\begin{Theorem}\label{3existence solution}
Assume that for all 
$(\theta,e)\in\Delta_{n}\times E^n$,
 the Brownian BSDE
\beq \nonumber 
Y_{t}^n(\theta,e) & = & \xi^n(\theta,e)+\int_{t}^Tf^n\Big(s,Y_{s}^n(\theta,e),Z_{s}^n(\theta,e),0,\theta,e\Big)ds\\
& & -\int_{t}^TZ_{s}^n(\theta,e)dW_{s}\;, \quad \theta_n\wedge T\leq t\leq T\;, \label{3existence BSDE n}
\enq
admits a solution $(Y^n(\theta,e),Z^n(\theta,e))\in \Sc^\infty_\F[\theta_n\wedge T,T]\times L^2_{\F}[\theta_n\wedge T,T]$, 
and that for each $k=0, \ldots,n-1$, the Brownian BSDE
\begin{equation} \label{3existence BSDE k}
\begin{aligned}
Y_{t}^k(\theta_{(k)}, e_{(k)}) ~ = & ~~ \xi^k(\theta_{(k)}, e_{(k)}) +\int_{t}^T f^k\Big(s,Y_{s}^k(\theta_{(k)}, e_{(k)}),Z_{s}^k(\theta_{(k)}, e_{(k)}), \\
&~~Y_{s}^{k+1}(\theta_{(k)},s, e_{(k)},.)-Y^k_{s}(\theta_{(k)}, e_{(k)}), \theta_{(k)}, e_{(k)}\Big)ds\\
 & ~ -\int_{t}^TZ_{s}^k(\theta_{(k)}, e_{(k)})dW_{s}\;,\quad \theta_k\wedge T\leq t\leq T\;,
\end{aligned}
\end{equation}
admits a solution $\big(Y^k(\theta_{(k)}, e_{(k)}), Z^k(\theta_{(k)}, e_{(k)})\big)\in \Sc^\infty_\F[\theta_k\wedge T,T]\times L^2_{\F}[\theta_k\wedge T,T]$. 
Assume moreover that each $Y^k$ (resp. $Z^k$) is $\PMc({\F})\otimes\Bc(\Delta_{k})\otimes\Bc(E^k)$-measurable (resp. $\Pc({\F})\otimes\Bc(\Delta_{k})\otimes\Bc(E^k)$-measurable).

If all these solutions  satisfy 
\beq\label{condSinfty}
\sup_{(k,\theta,e)
}{\|Y^k(\theta_{(k)},e_{(k)})\|}_{\Sc^\infty[\theta_k\wedge T,T]} & < & \infty \,,
\enq 
and 
\beqs
\E\Big[\int_{\Delta_{n}\times E^n}\Big(\int_{0}^{\theta_{1}\wedge T}|Z^0_{s}|^2ds+\sum_{k=1}^{n}\int_{\theta_{k}\wedge T}^{\theta_{k+1}\wedge T}|Z^k_{s}(\theta_{(k)},e_{(k)})|^2 ds \Big)\gamma_T(\theta,e)d\theta de \Big]  & < & \infty\;,
\enqs
then, under (HD), BSDE \reff{3BSDE jump} admits a solution $(Y, Z, U)\in \Sc^\infty_{\G}[0,T]\times L^2_\G[0,T]\times L^2(\mu)$  given by 
\begin{equation}\label{3expression solution} \left\{
\begin{aligned}
Y_{t}  ~& =~ Y^0_t\mathds{1}_{t< \tau_1}+\sum_{k=1}^{n}Y^k_t(\tau_{(k)}, \zeta_{(k)})\mathds{1}_{\tau_k\leq t< \tau_{k+1}}\;,\\
Z_{t} ~& =  ~Z^0_t\mathds{1}_{t\leq \tau_1}+\sum_{k=1}^{n}Z^k_t(\tau_{(k)}, \zeta_{(k)})\mathds{1}_{\tau_k< t\leq \tau_{k+1}}\;,\\
U_{t}(.)~ & = ~ U^0_t(.)\mathds{1}_{t\leq \tau_1}+\sum_{k=1}^{n-1}U^k_t(\tau_{(k)}, \zeta_{(k)},.)\mathds{1}_{\tau_k< t\leq \tau_{k+1}}\;,
\end{aligned}
\right.\end{equation}
with $U^k_{t}(\tau_{(k)}, \zeta_{(k)},.)=Y^{k+1}_{t}(\tau_{(k)}, t, \zeta_{(k)},.)-Y^{k}_{t}(\tau_{(k)}, \zeta_{(k)})$ for each $k = 0, \ldots, n-1$.    
\end{Theorem}

\vspace{2mm}

\noindent\textbf{Proof.} 
To alleviate notation, we shall often write $\xi^k$ and $f^k(t, y, z, u)$ instead of $\xi^k(\theta_{(k)}, e_{(k)})$ and $f^k(t, y, z, u, \theta_{(k)}, e_{(k)})$, and $Y^k_t(t, e)$ instead of $Y^k_t(\theta_{(k-1)}, t, e_{(k-1)}, e)$.

\noindent\textbf{Step 1:} We prove that for $t$ $\in$ $[0,T]$, $(Y, Z, U)$ defined by \reff{3expression solution} satisfies the equation  
\beq\label{relamontrer}
Y_t & = & \xi+\int_t^Tf(s,Y_s,Z_s,U_s)ds-\int_t^TZ_sdW_s-\int_t^T\int_{E}U_s(e)\mu(de,ds)\;.
\enq
We make an induction on the number $k$ of jumps in $(t,T]$.

\noindent$\bullet$ Suppose that $k$ $=$ $0$. 
We distinguish  two cases. \\
\textbf{Case 1:} there are $n$ jumps before $t$. We then have $\tau_n \leq t$ and from \reff{3expression solution} we get $Y_t=Y^n_t$. Using BSDE \reff{3existence BSDE n}, we can see that
\beqs
Y_t& = &\xi^n+\int_t^Tf^n(s,Y^n_s,Z^n_s,0)ds-\int_t^TZ^n_sdW_s\;.
\enqs
Since $\tau_n \leq T$, we have $\xi^n = \xi$ from \reff{3decomposition xi}. In the same way, we have $Y_s=Y^n_s$, $Z_s = Z^n_s$ and $U_s=0$ for all $s \in (t, T]$ from \reff{3expression solution}. Using \reff{3decomposition f}, we also get $f^n(s,Y^n_s,Z^n_s,0)=f(s,Y_s,Z_s,U_s)$ for all $s \in (t, T]$. Moreover, since the predictable processes $Z\mathds{1}_{\tau_n< .}$ and $Z^n\mathds{1}_{\tau_n< .}$ are indistinguishable on $\{\tau_n\leq t\}$, we have from Theorem 12.23 of \cite{HeW},  $\int_t^TZ_sdW_s = \int_t^TZ_s^ndW_s $ on $\{\tau_n \leq t\}$. 
Hence, we get
\beqs
Y_t & = & \xi+\int_t^Tf(s,Y_s,Z_s,U_s)ds-\int_t^TZ_sdW_s-\int_t^T\int_{E}U_s(e)\mu(de,ds)\;,
\enqs
 on $\{\tau_n \leq t\}$.\\
\noindent \textbf{Case 2:} there are $i$ jumps before $t$ with $i<n$ hence $Y_t = Y^i_t$. 
Since there is no jump after $t$, we have $Y_s=Y^i_s$, $Z_s = Z^i_s$, $U^i_{s}(.)=Y^{i+1}_{s}(s,.)-Y^{i}_{s}$, $\xi=\xi^i$ and $f^i(s,Y^i_s,Z^i_s,U^i_s)=f(s,Y_s,Z_s,U_s)$ for all $s \in (t,T]$, 
and  $\int_t^T\int_{E}U_s(e)\mu(de,ds)=0$.  Since the predictable processes $Z\mathds{1}_{\tau_i< .\leq \tau_{i+1}}$ and $Z^i\mathds{1}_{\tau_i< .\leq \tau_{i+1}}$ are indistinguishable on $\{\tau_i\leq t\}\cap\{T<\tau_{i+1}\}$, we have from Theorem 12.23 of \cite{HeW},  $\int_t^TZ_sdW_s = \int_t^TZ_s^idW_s $ on $\{\tau_i \leq t \}\cap\{T <\tau_{i+1} \}$. Combining these equalities with \reff{3existence BSDE k}, we get 
\beqs
Y_t&=&\xi+\int_t^Tf(s,Y_s,Z_s,U_s)ds-\int_t^TZ_sdW_s-\int_t^T\int_{E}U_s(e)\mu(de,ds)\;,
\enqs
on $\{\tau_i \leq t \}\cap\{T <\tau_{i+1} \}$.\\

\noindent $\bullet$ Suppose equation \reff{relamontrer} holds true when there are $k$ jumps in $(t,T]$, and consider the case where there are $k+1$ jumps in 
$(t,T]$. 

Denote by $i$ the number of jumps in $[0,t]$ hence $Y_t = Y^i_t$. Then, we have $Z_s = Z^i_s$, $U^i_{s}(.)=Y^{i+1}_{s}(s,.)-Y^{i}_{s}$ for all $s \in (t,\tau_{i+1}]$, and $Y_s=Y^i_s$ and $f(s,Y_s,Z_s,U_s)$ $=$ $f^i(s,Y^i_s,Z^i_s,U^i_s)$ for all $s \in (t,\tau_{i+1})$. Using \reff{3existence BSDE k}, we have
\beqs\nonumber
Y_{t}&=&Y^{i}_{\tau_{i+1}}+\int_t^{\tau_{i+1}}f(s,Y_s,Z_s,U_s)ds
-\int_t^{\tau_{i+1}}Z_s^idW_s\\\nonumber
 & = & Y^{i+1}_{\tau_{i+1}} +\int_t^{\tau_{i+1}}f(s,Y_s,Z_s,U_s)ds
-\int_t^{T}Z_s^i\mathds{1}_{\tau_i< s\leq \tau_{i+1}}dW_s \\
&&-\int_t^{\tau_{i+1}}\int_EU_s(e)\mu(de,ds)\;.
\enqs
Since the predictable processes $Z\mathds{1}_{\tau_i<.\leq\tau_{i+1}}$ and $Z^i\mathds{1}_{\tau_i<.\leq \tau_{i+1}}$ are indistinguishable on $\{\tau_i \leq t < \tau_{i+1}\}\cap\{\tau_{i+k+1}\leq T < \tau_{i+k+2}\}$, we get from Theorem 12.23 of \cite{HeW}, that $\int_t^{T}Z_s^i\mathds{1}_{\tau_i< s \leq \tau_{i+1}}dW_s=\int_t^{T}Z_s^i\mathds{1}_{\tau_i< s \leq \tau_{i+1}}dW_s$. Therefore, we get
\beq\label{eqttaui+1}
Y_{t} & = & Y^{i+1}_{\tau_{i+1}} +\int_t^{\tau_{i+1}}f(s,Y_s,Z_s,U_s)ds
-\int_t^{\tau_{i+1}}Z_sdW_s 
-\int_t^{\tau_{i+1}}\int_EU_s(e)\mu(de,ds)\;,\qquad \label{eqttaui+1}
\enq
on $\{\tau_i \leq t < \tau_{i+1}\}\cap\{\tau_{i+k+1}\leq T< \tau_{i+k+2}\}$.
Using the induction assumption on $(\tau_{i+1},T]$, we have 
\beqs
Y_{r}\mathds{1}_A(r) & = & \Big(\xi+\int_{r}^T\hspace{-2mm}f(s,Y_s,Z_s,U_s)ds
-\int_{r}^T\hspace{-2mm}Z_sdW_s-\int_{r}^T\hspace{-1mm}\int_E\hspace{-1mm}U_s(e)\mu(de,ds)\Big)\mathds{1}_A(r)\;,\enqs
for all $r\in[0,T]$, where 
\beqs
A= \Big\{(\omega,s)\in\Omega\times[0,T]~:~\tau_{i+1}(\omega)\leq s < \tau_{i+2}(\omega) ~ \mbox{ and } \tau_{i+k+1}(\omega)\leq T < \tau_{i+k+2}(\omega) \Big\}\;.
\enqs
Thus, the processes $Y_{}\mathds{1}_A(.)$ and $\Big(\xi+\int_{.}^T\hspace{-2mm}f(s,Y_s,Z_s,U_s)ds
-\int_{.}^T\hspace{-2mm}Z_sdW_s-\int_{.}^T\hspace{-2mm}\int_E\hspace{-0,5mm}U_s(e)\mu(de,ds)\Big)\mathds{1}_A(.)$
are indistinguishable since they are c\`ad-l\`ag modifications of the other. In particular they coincide at the stopping time $\tau_{i+1}$ and we get from the definition of $Y$
\beq
Y^{}_{\tau_{i+1}} ~=~Y^{i+1}_{\tau_{i+1}} & = & \xi+\int_{\tau_{i+1}}^T\hspace{-2mm}f(s,Y_s,Z_s,U_s)ds
-\int_{\tau_{i+1}}^T\hspace{-2mm}Z_sdW_s\nonumber\\
 & & 
-\int_{\tau_{i+1}}^T\hspace{-1mm}\int_E\hspace{-1mm}U_s(e)\mu(de,ds)\;.\quad \label{eqtuai+1T}
\enq
Combining \reff{eqttaui+1} and \reff{eqtuai+1T}, we get \reff{relamontrer}.
\vspace{2Mm}

\noindent\textbf{Step 2:} Notice that the process $Y$  (resp. $Z$, $U$) is $\Pc\Mc(\G)$ (resp. $\Pc(\G)$, $\Pc(\G)\otimes\Bc(E)$)-measurable since each $Y^k$ (resp. $Z^k$) is $\PMc({\F})\otimes\Bc(\Delta_{k})\otimes\Bc(E^k)$ (resp. $\Pc({\F})\otimes\Bc(\Delta_{k})\otimes\Bc(E^k)$)-measurable.

\vspace{2mm}

\noindent\textbf{Step 3:} We now prove that the solution satisfies the integrability conditions. Suppose that  the processes $Y^k$, $k=0,\ldots,n$, satisfy \reff{condSinfty}.
Define the constant $M$ by 
\beqs
M & := & \sup_{(k,\theta,e)
}{\|Y^k(\theta_{(k)},e_{(k)})\|}_{\Sc^\infty[\theta_k\wedge T,T]} \;,
\enqs
and consider the set $A\in\Fc_T\otimes\Bc(\Delta_n\cap[0,T]^n)\otimes\Bc(E^n)$ defined by
\beqs
A & := & \Big\{(\omega,\theta,e)\in\Omega\times(\Delta_n\cap[0,T]^n)\times E^n~:~\max_{0\leq k\leq n}\sup_{t\in[\theta_k,T]}|Y^k_t(\theta_{(k)},e_{(k)})|~\leq~M\Big\}\;.
\enqs
Then, we have $\P(\tilde\Omega)=1$, where 
\beqs
\tilde \Omega & = & \Big\{\omega\in\Omega~:~\big(\omega,\tau(\omega),\zeta(\omega)\big)\in A \Big\} \;.
\enqs
Indeed, we have from the density assumption (HD)
\beq\nonumber
\P(\tilde\Omega^c) & = & \E\Big[\mathds{1}_{A^c}\big(\omega,\tau(\omega),\zeta(\omega)\big)\Big]~=~ \E\Big[\E\Big[\mathds{1}_{A^c}\big(\omega,\tau(\omega),\zeta(\omega)\big)\Big|\Fc_T\Big]\Big]\\
 & = & \int_{(\Delta_n\cap[0,T]^n)\times E^n}\E\Big[\mathds{1}_{A^c}\big(\omega,\theta,e\big)\gamma_T(\theta,e)\Big]d\theta de\;.\label{tildeomegapleine proba}
\enq
From the definition of $M$ and $A$, we have 
\beqs
\mathds{1}_{A^c}\big(.,\theta,e\big)\gamma_T(\theta,e) & = & 0 \;,\quad \P-a.s. \;,
\enqs
for all $(\theta,e)\in (\Delta_n\cap[0,T]^n)\times E^n$.
Therefore, we get from \reff{tildeomegapleine proba}, $\P(\tilde\Omega^c)=0$. 
Then, by definition of $Y$, we have
\beqs
|Y_t| & \leq &  \big|Y^0_t\big|\mathds{1}_{t< \tau_1}+\sum_{k=1}^{n}\big|Y^k_t(\tau_{(k)}, \zeta_{(k)})\big|\mathds{1}_{\tau_k\leq t}\;.
\enqs
Since $\P(\tilde \Omega)=1$,  we have 
\beq\label{maj mominfty Y}
\big|Y^k_t(\tau_{(k)}, \zeta_{(k)})\big|\mathds{1}_{\tau_k\leq t} & \leq  & M\;, \quad 0\leq k\leq n\;,\quad \P-a.s. 
\enq
Therefore, we get from \reff{maj mominfty Y} 
\beqs
|Y_t| & \leq & (n+1)M\;, \quad \P-a.s. \;,
\enqs
for all $t\in[0,T]$. Since $Y$ is c\`ad-l\`ag, we get 
\beqs
\|Y\|_{\Sc^\infty[0,T]} & \leq & (n+1)M\;. 
\enqs
In the same way, using (HD) and the tower property of conditional expectation, we get
\beqs
\E\Big[ \int_0^T |Z_s|^2ds \Big] = \E\Big[ \int_{\Delta_{n}\times E^n}\Big(\int_{0}^{\theta_{1}\wedge T}|Z^0_{s}|^2ds+\sum_{k=1}^{n}\int_{\theta_{k}\wedge T}^{\theta_{k+1}\wedge T}|Z^k_{s}(\theta_{(k)},e_{(k)})|^2 ds \Big)\gamma_T(\theta,e)d\theta de\Big] \;. 
\enqs
Thus, $Z\in L^2_\G[0,T]$ since the processes $Z^k$, $k= 0,\ldots,n$, satisfy
\beqs
\E\Big[\int_{\Delta_{n}\times E^n}\Big(\int_{0}^{\theta_{1}\wedge T}|Z^0_{s}|^2ds+\sum_{k=1}^{n}\int_{\theta_{k}\wedge T}^{\theta_{k+1}\wedge T}|Z^k_{s}(\theta_{(k)},e_{(k)})|^2 ds \Big)\gamma_T(\theta,e)d\theta de\Big]  & < & \infty\;.
\enqs
Finally, we check that $U$ $\in$ $L^2(\mu)$. Using (HD), we have 
\beqs
\|U\|_{L^2(\mu)}^2 & = & \sum_{k=1}^n\int_{\Delta_n\times E^n} \E\Big[|Y^k_{\theta_k}(\theta_{(k)},e_{(k)})-Y^{k-1}_{\theta_k}(\theta_{(k-1)},e_{(k-1)})|^2\gamma_T(\theta,e)\Big]d\theta de\\
 & \leq & 2 \sum_{k=1}^n
 \Big(\|Y^k(\theta_{(k)},e_{(k)})\|_{\Sc^\infty[\theta_k\wedge T,T]}^2+\|Y^{k-1}(\theta_{(k-1)},e_{(k-1)})\|_{\Sc^\infty[\theta_{k-1}\wedge T,T]}^2\Big)
 \\
 & < & \infty\;.
\enqs
Hence,  $U$ $\in$ $L^2(\mu)$.
\ep

\begin{Remark}{\rm From the construction of the solution of BSDE \reff{3BSDE jump}, the jump component $U$ is bounded in the following sense
\beqs
\sup_{e\in E} {\|U(e)\|}_{\Sc^\infty[0,T]} & < & \infty\;.
\enqs
In particular, the random variable $\esssup_{(t,e)\in[0,T]\times E}|U_t(e)|$ is bounded.    }
\end{Remark}
\subsection{Application to quadratic BSDEs with jumps}

We suppose that the random variable $\xi$  and the generator $f$ satisfy the following conditions:
\begin{enumerate}[(HEQ1)]
\item The random variable $\xi$ is bounded: there exists a positive constant $C$ such that 
\beqs
|\xi| & \leq & C\;, \quad\P-a.s. 
\enqs 
\item The generator $f$ is quadratic in $z$: there exists a constant  $C$ such that 
\beqs
|f(t,y,z,u)| & \leq & C \Big(1 + |y| + |z|^2 + \int_E |u(e)| \lambda_t(e) de \Big)\;,
\enqs
for all $(t,y,z,u) \in [0,T]\times\R\times\R^d \times Bor(E,\R)$.
\item For any $R>0$, there exists a function $mc^f_R$ such that $\lim_{\eps\rightarrow0}mc^f_R(\eps)=0$ and
\beqs
|f_t(y,z,(u(e)-y)_{e\in E})-f_t(y',z',(u(e)-y)_{e\in E})| & \leq & mc^f_R(\eps)\;,
\enqs 
for all $(t,y,y',z,z',u)\in [0,T]\times[\R]^2\times[\R^d]^2 \times Bor(E,\R)$ s.t. $|y|,|z|,|y'|,|z'|\leq R$ and $|y-y'|+|z-z'|\leq \eps$.
\end{enumerate}
\begin{Proposition}\label{exquad} 
Under  (HD), (HBI), (HEQ1), (HEQ2) and  (HEQ3), BSDE \reff{3BSDE jump} admits a solution in $\Sc^\infty_\G[0,T] \times L^2_\G[0,T] \times L^2(\mu)$.
\end{Proposition}

\noindent\textbf{Proof.}
\textbf{Step 1.} Since $\xi$ is a bounded random variable, we can choose $\xi^k$ bounded for each $k=0,\ldots,n$. Indeed, let $C$ be a positive constant such that $|\xi|\leq C,~\P-a.s.$, then, we have
\beqs
\xi & = & \sum_{k=0}^{n}\tilde \xi^k(\tau_{1}, \ldots, \tau_{k},  \zeta_1, \ldots,  \zeta_k)\mathds{1}_{\tau_k\leq T< \tau_{k+1}}\;,
\enqs
with $\tilde \xi ^k(\tau_{1}, \ldots, \tau_{k},  \zeta_1, \ldots,  \zeta_k)=(\xi^k(\tau_{1}, \ldots, \tau_{k},  \zeta_1, \ldots,  \zeta_k)\wedge C)\vee (-C)$, for each $k=0, \ldots, n$.

\vspace{2mm}

\noindent\textbf{Step 2.} Since $f$ is quadratic in $z$, it is possible to choose the functions $f^k$, $k =0, \ldots, n$, quadratic in $z$. Indeed, if $C$ is a positive constant such that  $|f(t,y,z,u)| \leq C(1 + |y| + |z|^2 + \int_E |u(e)| \lambda_t(e) de )$, for all $(t,y,z,u) \in [0,T]\times\R\times\R^d \times Bor(E,\R), ~\P-a.s.$ and $f$ has the following decomposition
\beqs
f(t,y,z,u) & = &\sum_{k=0}^{n}f^k(t, y, z, u, \tau_{(k)}, \zeta_{(k)})
\mathds{1}_{\tau_k\leq t< \tau_{k+1}}\;,
\enqs
then, $f$ satisfies the same decomposition with $\tilde f^k$ instead of $f^k$  where
\beqs
\tilde f ^k(t,y,z,u,\theta_{(k)},e_{(k)}) & = & f^k(t, y, z, u,\theta_{(k)},e_{(k)})\wedge \Big(C \Big(1 + |y| + |z|^2 + \int_E |u(e)| \lambda_t(e) de \Big)\Big) \\
& & \vee \Big(-C \Big(1 + |y| + |z|^2 + \int_E |u(e)| \lambda_t(e) de\Big) \Big) \;, 
\enqs
for all $(t,y,z,u)\in[0,T]\times\R\times\R^d \times Bor(E,\R)$ and $(\theta,e)\in\Delta_n\times E^n$.

\vspace{2mm}

\noindent\textbf{Step 3.}
We now prove by a backward induction that there exists for each $k = 0, \ldots, n-1$ (resp. $k=n$), a solution $(Y^k, Z^k)$ to BSDE \reff{3existence BSDE k} (resp. \reff{3existence BSDE n}) s.t. $Y^k$ is a $\Pc\Mc(\F)\otimes\Bc(\Delta_k)\otimes\Bc(E^k)$-measurable process and $Z^k$ is a $\Pc(\F)\otimes\Bc(\Delta_k)\otimes\Bc(E^k)$-measurable process, and 
\beqs
\sup_{(\theta_{(k)}, e_{(k)}) \in  \Delta_k \times E^k}{\|Y^k(\theta_{(k)}, e_{(k)})\|}_{\Sc^\infty[\theta_k\wedge T, T]} + {\|Z^k(\theta_{(k)}, e_{(k)})\|}_{L^2[\theta_k\wedge T, T]} & < & \infty \;.
\enqs

\vspace{1mm}

\ni $\bullet$ Choosing $\xi^n(\theta_{(n)}, e_{(n)})$ bounded as in Step 1, we get from (HEQ3) and Proposition \ref{reg-decomp} and Theorem 2.3 of \cite{kob00} the existence of a solution $(Y^n(\theta_{(n)}, e_{(n)}),Z^n(\theta_{(n)}, e_{(n)}))$ to BSDE \reff{3existence BSDE n}. 

We now check that we can choose  $Y^n$ (resp. $Z^n$) as a $\Pc\Mc(\F)\otimes\Bc(\Delta_n)\otimes\Bc(E^n)$ (resp. $\Pc(\F)\otimes\Bc(\Delta_n)\otimes\Bc(E^n)$)-measurable process. Indeed, we know (see \cite{kob00}) that we can construct the solution $(Y^n, Z^n)$ as limit of solutions to Lipschitz BSDEs. From Proposition \ref{mesLipsol}, we then get a $\Pc(\F)\otimes\Bc(\Delta_n)\otimes\Bc(E^n)$-measurable solution as limit of $\Pc(\F)\otimes\Bc(\Delta_n)\otimes\Bc(E^n)$-measurable processes.  Hence, $Y^n$ (resp. $Z^n$) is a $\Pc\Mc(\F)\otimes\Bc(\Delta_n)\otimes\Bc(E^n)$ (resp. $\Pc(\F)\otimes\Bc(\Delta_n)\otimes\Bc(E^n)$)-measurable process. 
Applying Proposition 2.1 of \cite{kob00} to $(Y^n, Z^n)$, we get from (HEQ1) and (HEQ2)
\beqs
\sup_{(\theta,e)\in \Delta_n \times E^n}{\|Y^n(\theta_{(n)},e_{(n)})\|}_{\Sc^\infty[\theta_n \wedge T, T]} + {\|Z^n(\theta_{(n)}, e_{(n)})\|}_{L^2[\theta_n \wedge T, T]} & < & \infty \;.
\enqs

\vspace{2mm}

\ni $\bullet$ Fix $k \leq n-1$ and suppose that the result holds true for $k+1$: there exists $(Y^{k+1}, Z^{k+1})$ such that 
\beqs
\sup_{(\theta_{(k+1)}, e_{(k+1)}) \in \Delta_{k+1} \times E^{k+1}}\Big\{{\|Y^{k+1}(\theta_{(k+1)}, e_{(k+1)})\|}_{\Sc^\infty[\theta_{k+1}\wedge T, T]}\quad  & & \\
+ {\|Z^{k+1}(\theta_{(k+1)}, e_{(k+1)})\|}_{L^2[\theta_{k+1}\wedge T, T]}\Big\} & < & \infty \;.
\enqs
 Then, using  (HBI), there exists a constant $C > 0$ such that 
\beqs
\Big|f^k\Big(s, y , z ,Y_{s}^{k+1}(\theta_{(k)},s, e_{(k)},.) - y), \theta_{(k)}, e_{(k)}\Big)\Big|
\leq C ( 1 + |y| + |z|^2 ) \;.
\enqs
Choosing $\xi^k(\theta_{(k)}, e_{(k)})$ bounded as in Step 1, we get from  (HEQ3) and Proposition \reff{reg-decomp} and Theorem 2.3 of \cite{kob00} the existence of a solution $(Y^k(\theta_{(k)}, e_{(k)}), Z^k(\theta_{(k)}, e_{(k)}))$. 

As for $k=n$, we can choose $Y^k$ (resp. $Z^k$) as a $\Pc\Mc(\F)\otimes\Bc(\Delta_k)\otimes\Bc(E^k)$ (resp. $\Pc(\F)\otimes\Bc(\Delta_k)\otimes\Bc(E^k)$)-measurable process. 

Applying Proposition 2.1 of \cite{kob00} to $(Y^k(\theta_{(k)}, e_{(k)}), Z^k(\theta_{(k)}, e_{(k)}))$, we get from (HEQ1) and (HEQ2)
\beqs
\sup_{(\theta_{(k)}, e_{(k)}) \in \Delta_k \times E^k}{\|Y^k(\theta_{(k)}, e_{(k)})\|}_{\Sc^\infty[\theta_k\wedge T, T]} + {\|Z^k(\theta_{(k)}, e_{(k)})\|}_{L^2[\theta_k\wedge T, T]} & < & \infty \;.
\enqs

\ni\textbf{Step 4.} From Step 3, we can apply Theorem \ref{3existence solution}. We then get the existence of a solution to BSDE \reff{3BSDE jump}.
\ep
\begin{Remark}
{ \rm  Our existence result is given for bounded terminal condition. It is based on the result of Kobylanski for quadratic Brownian BSDEs in \cite{kob00}. We notice that existence results for quadratic BSDEs with unbounded terminal conditions have recently been proved in  Briand and Hu \cite{bh06} and  Delbaen \textit{ et al.} \cite{delhuri11}. These works provide existence results for solutions of Brownian quadratic BSDEs with exponentially integrable terminal conditions and generators and conclude that the solution  $Y$ satisfies an exponential integrability condition.     

%
Here, we cannot use these results in our approach. Indeed, consider the case of a single jump with the generator $f(t,y,z,u)=|z|^2+|u|$. The associated decomposed BSDE at rank $0$ is given by 
\beqs
Y^0_t & = & \xi^0 + \int_t^T \big[|Z^0_s|^2 + |Y^1_s(s)-Y^0_s| \big]ds-\int_t^TZ^0_sdW_s\;,\quad t\in[0,T]\;. 
\enqs
Then to apply the results from \cite{bh06} or  \cite{delhuri11},  we require that the process $(Y^1_s(s))_{s\in[0,T]}$ satisfies some exponential integrability condition. However, at rank 1, the decomposed BSDE is given by
\beqs
Y^1_t(\theta) & = & \xi^1(\theta) + \int_t^T \big|Z^1_s(\theta)\big|^2ds-\int_t^TZ^1_s(\theta)dW_s\;,\quad t\in[\theta,T]\;, \quad\theta\in [0,T] \;,
\enqs
and since $\xi^1$ satisfies an exponential integrability condition by assumption
 we know that $Y^1(\theta)$ satisfies an exponential integrability condition for any $\theta\in [0,T]$, but we have no information about the process $(Y^1_s(s))_{s \in [0,T]}$. The difficulty here lies in understanding the behavior of the ``sectioned" process $\{Y^1_s(\theta)~:~s=\theta\}$ and its study is left for further research. } 
\end{Remark}

\subsection{Application to the pricing of a European option in a market with a jump}
In this example, we assume that $W$ is one dimensional ($d=1$) and there is a single random time $\tau$ representing the time of occurrence of a shock in the prices on the market. We denote by $H$ the associated pure jump process:
\beqs
H_t & = & \mathds{1}_{\tau \leq t}\;,\quad 0\leq t \leq T\;.
\enqs
We consider a financial market  which consists of 
\begin{itemize}
\item a non-risky asset $S^0$, whose strictly positive price process is defined by
\beqs
dS^0_t &  =  & r_t S^0_t dt\;,\quad 0\leq t \leq T\;, \quad  S^0_0=1\;,
\enqs
with $r_t  \geq 0$, for all $t\in[0,T]$,

\item  two risky assets with respective price processes $S^1$ and $S^2$ defined by
\beqs
dS^1_t & = & S^1_{t^-} (b_t dt + \sigma_t dW_t + \beta dH_t)\,, \quad 0\leq t\leq T\;,\quad  S^1_0 = s^1_0\;,
\enqs
and 
\beqs
dS^2_t & = & S^2_t (\bar{b}_t dt + \bar{\sigma}_t dW_t)\,, \quad 0\leq t\leq T\;,\quad   S^2_0 = s^2_0\,,
\enqs
with $\sigma_t > 0$ and $\bar{\sigma}_t >0$, and $\beta > -1$ (to ensure that the price process $S^1$ always remains strictly positive).
\end{itemize}
 We make the following assumption which ensures the existence of the processes $S^0$, $S^1$, and $S^2$:

\vspace{2mm}

\ni(HB) The coefficients $r$, $b$, $\bar b$, $\sigma$, $\bar \sigma$, $\frac{1}{\sigma}$ and $\frac{1}{\bar \sigma}$ are bounded: there exists a constant $C$ s.t.
\beqs
|r_t| + |b_t| + |\bar b_t| + |\sigma_t| + |\bar\sigma_t| + \Big|\frac{1}{\sigma_t}\Big| + \Big|\frac{1}{\bar \sigma_t}\Big| & \leq & C\;, ~0 \leq t \leq T\;, ~\P-a.s.
\enqs

\vspace{2mm}

We assume that the coefficients $r$, $b$, $\bar b$, $\sigma$ and $\bar \sigma$ have the following forms 
\[\left\{ \begin{aligned}
r_t = r^0 \mathds{1}_{t < \tau} + r^1(\tau) \mathds{1}_{t \geq \tau}\,,\\
b_t = b^0 \mathds{1}_{ t < \tau} + b^1(\tau)\mathds{1}_{ t \geq \tau}\,,\\
\bar{b}_t = \bar{b}^0 \mathds{1}_{ t < \tau} + \bar{b}^1(\tau)\mathds{1}_{ t \geq \tau}\,,\\
\sigma_t = \sigma^0 \mathds{1}_{ t < \tau} + \sigma^1(\tau)\mathds{1}_{ t \geq \tau}\,,\\
\bar{\sigma}_t = \bar{\sigma}^0 \mathds{1}_{ t < \tau} + \bar{\sigma}^1(\tau)\mathds{1}_{ t \geq \tau}\,,
\end{aligned}\right.\]
for all $t\geq 0$. 

The aim of this subsection is to provide an explicit price for any bounded $\Gc_{T}$-measurable European option $\xi$ of the form
\beqs
\xi & = & \xi^0\mathds{1}_{ T<\tau}+\xi^1(\tau)\mathds{1}_{ \tau\leq T}\;,
\enqs
where $\xi^0$ is $\Fc_T$-measurable and $\xi^1$ is $\Fc_T\otimes\Bc(\R)$-measurable, 
 together with a replicating strategy  $\pi$ $=$ $(\pi^0, \pi^1, \pi^2)$ ($\pi^i_t$ corresponds to the number of share of $S^i$ held at time $t$). 
We assume that this market model is free of arbitrage oppotunitity (a necessary and sufficient condition to ensure it is e.g. given in Lemma 3.1.1 of \cite{gio10}).
%
%
%
%
%

The value of a contingent claim is then given by the initial amount of a replicating portfolio. 
Let $\pi$ $=$ $(\pi^0, \pi^1, \pi^2)$ be a $\Pc(\G)-$measurable self-financing strategy. The wealth process $Y$ associated with this strategy satisfies
\beq\label{3Ypi}
 Y_t & = & \pi^0_t S^0_t +\pi^1_t S^1_t +\pi^2_t S^2_t\,,\quad 0\leq t\leq T\;.
\enq
Since $\pi$ is a self-financing strategy, we have 
\beqs
dY_t & = & \pi^0_t dS^0_t + \pi^1_t dS^1_t + \pi^2_t dS^2_t\,,\quad 0\leq t\leq T\;.
\enqs 
Combining this last equation with \reff{3Ypi}, we get
\beq \nonumber
dY_t &  = &  \big( r_t Y_t +(b_t - r_t) \pi^1_t S^1_t + (\bar{b}_t -  r_t) \pi^2_t S^2_t  \big) dt\\
 & &  + \big( \pi^1_t \sigma_t S^1_{t} + \pi^2_t \bar{\sigma}_t S^2_t \big) dW_t +  \pi^1_t \beta S^1_{t^-} dH_t \,, \quad 0\leq t\leq T\;.\label{3equation richesse}
\enq
Define the predictable processes $Z$ and $U$ by
\beq\label{3Zpi}
Z_t = \pi^1_t \sigma_{t} S^1_{t} + \pi^2_t \bar{\sigma}_{t} S^2_t & \text{and} & U_t = \pi^1_t \beta S^1_{t^-}\;, \quad 0\leq t\leq T\,.
\enq
Then, \reff{3equation richesse} can be written under the form
\beqs
dY_t  &  = &  \Big[ r_t Y_t - \frac{r_t - \bar{b}_t}{ \bar{\sigma}_t} Z_t -\Big(\frac{r_t - b_t}{\beta}  - \frac{\sigma_t (r_t - \bar{b}_t)}{\beta\bar{\sigma}_t} \Big)U_t \Big] dt + Z_t dW_t +  U_tdH_t\;, \quad 0\leq t\leq T\;.
\enqs
Therefore, the problem of valuing and hedging of the contingent claim $\xi$ consists  in solving the following BSDE
\begin{equation}\label{3edsr pricing}
\left\{\begin{array}{rcl}
- dY_t & = &  \Big[ \frac{r_t - \bar{b}_t}{ \bar{\sigma}_t} Z_t +\Big(\frac{r_t - b_t}{\beta}  - \frac{\sigma_t (r_t - \bar{b}_t)}{\beta\bar{\sigma}_t} \Big)U_t - r_t Y_t\Big] dt \\
 & & - Z_t dW_t -  U_tdH_t\;,\quad 0\leq t\leq T\,,\\
Y_T & = &\xi\,.
\end{array}\right.
\end{equation}
The recursive system of Brownian BSDEs associated with  \reff{3edsr pricing} is then given by
\begin{equation}\label{3edsr apres saut}
\left\{\begin{array}{rcl}
- dY^1_t(\theta) & =  & \Big[ \frac{r^1(\theta) - \bar{b}^1(\theta)}{ \bar{\sigma}^1(\theta)} Z^1_t(\theta) - r^1(\theta) Y^1_t(\theta)\Big] dt - Z^1_t(\theta) dW_t \,,\quad \theta\leq t\leq T\;,\\
Y^1_T(\theta) & = & \xi^1 (\theta)\,,
\end{array}\right.
\end{equation}
and
\begin{equation}\label{3edsr avant saut}
\left\{\begin{array}{rcl}
- dY^0_t & = &  \Big[  \frac{r^0 - \bar{b}^0}{ \bar{\sigma}^0} Z_t + \Big(\frac{r^0 - b^0}{\beta}  - \frac{\sigma^0 (r^0 - \bar{b}^0)}{\beta\bar{\sigma}^0} \Big)(Y^1_t(t)-Y^0_t ) - r^0 Y^0_t\Big] dt\\
 & &  - Z_t dW_t\,,\qquad 0\leq t\leq T\;, \\
Y^0_T & = & \xi^0\,.
\end{array}\right.
\end{equation}

\begin{Proposition}
Under  (HD) and (HB), BSDE \reff{3edsr pricing} admits a solution in $\Sc_\G^\infty[0,T]\times L^2_{\G}[0,T]\times L^2(\mu)$.
\end{Proposition}
\ni\textbf{Proof.} Using the same argument as in Step 1 of the proof of  Proposition \ref{exquad}, we can assume w.l.o.g. that the coefficients of BSDEs \reff{3edsr apres saut} and \reff{3edsr avant saut} are bounded. Then, BSDE \reff{3edsr apres saut} is a linear BSDE with bounded coefficients and a bounded terminal condition. From Theorem 2.3 in \cite{kob00}, we get the existence of a solution $(Y^1(\theta),Z^1(\theta))$ in $\Sc^\infty_\F[\theta,T]\times L_\F^2[\theta,T]$ to  \reff{3edsr apres saut} for all $\theta\in[0,T]$. 
Moreover, from Proposition 2.1 in \cite{kob00}, we have
\beq\label{borne Y^1}
\sup_{\theta\in[0,T]}{ \|Y^1(\theta)\|}_{\Sc^\infty[\theta,T]} & < & \infty\;.
\enq
Applying Proposition \ref{mesLipsol} with $\Xc=[0,T]$ and $d\rho(\theta)=\gamma_0(\theta)d\theta$  we can choose the solution $(Y^1,Z^1)$ as a $\Pc(\F)\otimes\Bc([0,T])-$measurable process. 

Estimate \reff{borne Y^1} gives that BSDE \reff{3edsr avant saut} is also a linear BSDE with bounded coefficients. Applying Theorem 2.3  and Proposition 2.1 in \cite{kob00} as previously, we get the existence of a solution $(Y^0,Z^0)$ in $\Sc^\infty_\F[0,T]\times L_\F^2[0,T]$ to  \reff{3edsr avant saut}.  
Applying Theorem \ref{3existence solution}, we get the result.
\ep

\vspace{2mm}

Since BSDEs \reff{3edsr apres saut} and \reff{3edsr avant saut}  are linear, we have  explicit formulae for the solutions. 
For $Y^1(\theta)$, we get:
\beqs
Y^1_t(\theta) & = & \frac{1}{\Gamma^{1}_t(\theta)} \E\Big[ \xi^1(\theta) \Gamma^1_T(\theta) \Big| \Fc_t \Big]\,,\qquad \theta\leq t\leq T\;, 
\enqs 
with $\Gamma^1 (\theta)$ defined by
\beqs
\Gamma^1_t (\theta) & = & \exp\Big( \frac{r^1(\theta) - \bar{b}^1(\theta)}{ \bar{\sigma}^1(\theta)} W_t  - \frac{1}{2}  \Big|\frac{r^1(\theta) - \bar{b}^1(\theta)}{ \bar{\sigma}^1(\theta)} \Big|^2 t - r^1(\theta) t \Big)\,,\qquad \theta\leq t\leq T\,. 
\enqs
For $Y^0$, we get :
\beqs
Y^0_t & = & \frac{1}{{\Gamma}^{0}_t} \E \Big[ \xi^0 \Gamma^0_T + \int_t^T c_s \Gamma^0_s ds \Big| \Fc_t \Big]\,,\qquad 0\leq t\leq T\,,
\enqs
with $\Gamma^0$ defined by
\beqs
{\Gamma^0_t} & = & \exp \Big( \int_0^t d_s dW_s - \frac{1}{2} \int_0^t |d_s|^2 ds + \int_0^t a_s ds \Big)\,,\qquad 0\leq t\leq T\,,
\enqs
where the parameters $a$, $d$ and $c$ are given by 
\[\left\{\begin{aligned}
a_t & = -r^0 - \Big(\frac{r^0 - b^0}{\beta}  - \frac{\sigma^0 (r^0 - \bar{b}^0)}{\beta\bar{\sigma}^0} \Big)\,,\\
d_t & = \frac{r^0 - \bar{b}^0}{ \bar{\sigma}^0}\,, \\
c_t & = \Big(\frac{r^0 - b^0}{\beta}  - \frac{\sigma^0 (r^0 - \bar{b}^0)}{\bar{\beta\sigma}^0} \Big)Y^1_t(t)\,.
\end{aligned}\right.\]
The price at time $t$ of the European option $\xi$ is equal to $Y^0_t$ if $t<\tau$ and $Y^1_t(\tau)$ if $t\geq\tau$.  Once we know the processes $Y$ and $Z$, a hedging strategy $\pi=(\pi^0,\pi^1,\pi^2)$ is  given by   \reff{3Ypi} and \reff{3Zpi}.

\vspace{2mm}

Under no free lunch assumption,  all the hedging portfolios have the same value, which gives the uniqueness of the process $Y$.  This leads to the uniqueness issue for the whole solution $(Y,Z,U)$.

\section{Uniqueness}
\setcounter{equation}{0} \setcounter{Assumption}{0}
\setcounter{Theorem}{0} \setcounter{Proposition}{0}
\setcounter{Corollary}{0} \setcounter{Lemma}{0}
\setcounter{Definition}{0} \setcounter{Remark}{0}
In this section, we provide a uniqueness result based on a comparison theorem. 
We first provide a general comparison theorem which allows to compare solutions to the studied BSDEs  as soon as we can compare solutions to the associated system of recursive Brownian BSDEs. We then illustrate our general result with a concrete example in a convex framework.

\subsection{The general comparison theorem}
We consider two BSDEs with coefficients $(\underline f,\underline \xi)$ and $(\bar{f}, \bar{\xi})$ such that
\begin{itemize}
\item $\underline \xi$ (resp. $\bar \xi$) is a bounded $\Gc_T$-measurable random variable of the form
\beqs
\underline \xi & = &\sum_{k=0}^{n}\underline \xi^k(\tau_{(k)}, \zeta_{(k)})\mathds{1}_{\tau_k\leq T< \tau_{k+1}} \\
(\text{resp. }\bar \xi & = &\sum_{k=0}^{n}\bar \xi^k(\tau_{(k)}, \zeta_{(k)})\mathds{1}_{\tau_k\leq T< \tau_{k+1}}) \;,
\enqs
where $\underline \xi^0$ (resp. $\bar \xi^0$) is $\Fc_{T}$-measurable and $\underline \xi^k$ (resp. $\bar \xi^k$) is $\Fc_{T}\otimes\Bc(\Delta_{k})\otimes\Bc(E^{k})$-measurable for each $k=1,\ldots,n$, 
\item  $\underline f$ (resp. $\bar f$) is map from $[0,T]\times\Omega\times\R\times\R^d\times Bor(E,\R)$ to $\R$ which is  a $\Pc(\G)\otimes\Bc(\R)\otimes\Bc(\R^d)\otimes\Bc(Bor(E,\R))$-$\Bc(\R)$-measurable map. 
\end{itemize}

 We denote by $(\underline Y,\underline Z,\underline U)$ and $(\bar Y,\bar Z,\bar U)$ their respective solutions in $\Sc^\infty_{\G}[0,T]\times L^2_{\G}[0,T]\times L^2(\mu)$. 
We consider the decomposition 
 $(\underline Y^k)_{0\leq k\leq n}$ (resp. $(\bar Y^k)_{0\leq k\leq n}$, $(\underline Z^k)_{0\leq k\leq n}$, $(\bar Z^k)_{0\leq k\leq n}$, $(\underline U^k)_{0\leq k\leq n}$, $(\bar U^k)_{0\leq k\leq n}$ ) of $\underline Y$ (resp. $\bar Y$, $\underline Z$, $\bar Z$, $\underline U$,  $\bar U$) given by Lemma \ref{3lemme decomposition}.
For ease of notation, we shall write
$\underline F^k(t,y,z)$ and $\bar{F}^k(t,y,z)$ instead of $\underline f(t,y,z,\underline Y^{k+1}_t(\tau_{(k)},t,\zeta_{(k)},.)-y)$ and $\bar{f}(t,y,z,\bar Y^{k+1}_t(\tau_{(k)},t,\zeta_{(k)},.)-y)$ for each $k=0,\ldots,n-1$, and 
$\underline F^n(t,y,z)$ and $\bar{F}^n(t,y,z)$ instead of $\underline f(t,y,z,0)$ and $\bar{f}(t,y,z,0)$.

\ni We shall make, throughout the sequel, the standing assumption known as \textbf{(H)}-hypothesis:\\

(HC) Any $\F$-martingale remains a $\G$-martingale.

\vspace{3mm}

\begin{Remark}\label{rem W G-MG}
{\rm Since $W$ is an $\F-$Brownian motion, we get under (HC) that it remains a $\G-$Brownian motion. Indeed, using (HC), we have that $W$ is a $\G$-local martingale with quadratic variation $\langle W, W \rangle_t = t$. Applying L\'evy's characterization of Brownian motion (see e.g. Theorem 39 in \cite{pro05}), we obtain that $W$ remains a $\G$-Brownian motion. }
\end{Remark}



\begin{Definition}\label{defcompB}
{\rm We say that a generator $g:~\Omega\times[0,T]\times\R\times\R^d\rightarrow\R$ satisfies a comparison theorem for Brownian BSDEs if for any bounded $\G$-stopping times $\nu_2\geq \nu_1$, any generator $g':~\Omega\times[0,T]\times\R\times\R^d\rightarrow\R$ and any $\Gc_{\nu_2}$-measurable r.v.  $\zeta$ and $\zeta'$ such that $g\leq g'$ and $\zeta\leq \zeta'$ (resp. $g\geq g'$ and $\zeta\geq \zeta'$), we have $Y\leq Y'$ (resp. $Y\geq Y'$ ) on $[\nu_1,\nu_2]$. Here, $(Y,Z)$ and $(Y',Z')$ are solutions in $\Sc_\G^\infty[0,T]\times L^2_\G[0,T]$ to BSDEs with data $(\zeta,g)$ and $(\zeta',g')$:
\beqs
Y_t & = & \zeta+\int_t^{\nu_2}g(s,Y_s,Z_s)ds-\int_t^{\nu_2}Z_sdW_s\;,\quad \nu_1\leq t\leq \nu_2\;,
\enqs
and
\beqs
Y'_t & = & \zeta'+\int_t^{\nu_2}g'(s,Y'_s,Z'_s)ds-\int_t^{\nu_2}Z'_sdW_s\;,\quad \nu_1\leq t\leq \nu_2\;.
\enqs
 } 
\end{Definition}

\ni We can state the general  comparison theorem.
\begin{Theorem}\label{3theoreme comparaison lipschitz}
Suppose that $\underline \xi \leq \bar \xi,~\P$-a.s. 
Suppose moreover that for each $k=0,\ldots,n$ 
\beqs
\underline{F}^k(t,y,z) & \leq &\bar F^k(t,y,z),\quad  \forall (t,y,z)\in [0,T] \times \R \times \R^d,\quad\P-a.s.\;,
\enqs
and the generators $\bar{F}^k$ or $\underline{F}^k$ satisfy a comparison theorem for Brownian BSDEs. Then,  if $\bar U_{t}$ $=$ $\underline U_{t}$ $=$ $0$ for $t$ $>$ $\tau_{n}$, we have under (HD) and (HC) 
 \beqs
\underline{Y}_t  & \leq & \bar Y_t\;,\quad  0 \leq t \leq T\;,~~ ~\P-a.s.
 \enqs
\end{Theorem}

\noindent\textbf{Proof.} 
The proof is performed in four steps. We first identify the BSDEs of which the terms appearing in the decomposition of $\bar Y$ and $\underline Y$ are solutions in the filtration $\G$. We then modify $\bar Y ^k$ and $\underline Y ^k$ outside of $[\tau_k,\tau_{k+1})$ to get c\`ad-l\`ag processes for each $k= 0, \ldots, n$.  We then compare the modified processes by killing their jumps. Finally, we retrieve a comparison for the initial processes since the modification has happened outside of $[\tau_k,\tau_{k+1})$ (where they coincide with $\bar Y $ and $\underline Y$).

\vspace{2mm}

\noindent\textbf{Step 1.} Since $(\bar Y, \bar Z, \bar U)$ (resp. $(\underline Y, \underline Z, \underline U)$) is solution to the BSDE with parameters $(\bar\xi, \bar f)$ (resp. $(\underline \xi, \underline f)$), we obtain from the decomposition in the filtration $\F$  and Theorem 12.23 in \cite{HeW} that $(\bar{Y}^n, \bar{Z}^n)$ (resp. $(\underline Y^n, \underline Z^n)$) is solution to 
\beq \nonumber
\bar Y^n_{t}(\tau_{(n)},\zeta_{(n)}) & = & \bar  \xi +\int_{t}^T \bar{F}^n\Big(s,\bar Y^n_{s}(\tau_{(n)},\zeta_{(n)}),\bar Z^n_{s}(\tau_{(n)},\zeta_{(n)})\Big)ds\qquad\\\label{3EDSRBn+}
& &  -\int_{t}^{T}\bar Z^n_{s}(\tau_{(n)},\zeta_{(n)})dW_{s}\;,\quad \tau_{n} \wedge T \leq t \leq T\;,\\\nonumber
\\\nonumber
({\rm resp. }\quad \underline  Y^n_{t}(\tau_{(n)},\zeta_{(n)}) & = & \underline \xi+\int_{t}^{T}\underline F^n\Big(s,\underline Y^n_{s}(\tau_{(n)},\zeta_{(n)}),\underline Z^n_{s}(\tau_{(n)},\zeta_{(n)})\Big)ds\qquad\nonumber\\
& &  -\int_{t}^{T}\underline Z^n_{s}(\tau_{(n)},\zeta_{(n)})dW_{s}\;,\quad \tau_{n} \wedge T \leq t \leq T\;)\label{3EDSRBn-}
\enq
and $(\bar Y^k,\bar Z^k)$ (resp. $(\underline Y^k,\underline Z^k)$) is solution to
\beq\nonumber
\bar Y^k_{t}(\tau_{(k)},\zeta_{(k)}) & = & \big[ \bar  Y^{k+1}_{\tau_{k+1}}(\tau_{(k+1)},\zeta_{(k+1)}) - \bar U_{\tau_{k+1}}(\zeta_{k+1})\big] \mathds{1}_{\tau_{k+1} \leq T} + \bar \xi  \mathds{1}_{\tau_{k+1} > T} \\\nonumber
 & & +\int_{t}^{\tau_{k+1} \wedge T}\bar F^k\Big(s,\bar Y^k_{s}(\tau_{(k)},\zeta_{(k)}),\bar Z^k_{s}(\tau_{(k)},\zeta_{(k)})\Big)ds\qquad\\\label{3EDSRBk+}
& &  -\int_{t}^{\tau_{k+1}\wedge T}\bar Z^k_{s}(\tau_{(k)},\zeta_{(k)})dW_{s}\;,\quad \tau_{k} \wedge T \leq t < \tau_{k+1}\wedge T\;,\\\nonumber
\\\nonumber
({\rm resp. }\quad \underline  Y^k_{t}(\tau_{(k)},\zeta_{(k)}) & = & \big[\underline Y^{k+1}_{\tau_{k+1}}(\tau_{(k+1)},\zeta_{(k+1)}) - \underline U_{\tau_{k+1}}(\zeta_{k+1})\big] \mathds{1}_{\tau_{k+1} \leq T} + \underline \xi  \mathds{1}_{\tau_{k+1} > T}\\
 & & +\int_{t}^{\tau_{k+1}\wedge T} \underline F^k\Big(s,\underline Y^k_{s}(\tau_{(k)},\zeta_{(k)}),\underline Z^k_{s}(\tau_{(k)},\zeta_{(k)})\Big)ds\qquad\nonumber\\
& &  -\int_{t}^{\tau_{k+1}\wedge T}\underline Z^k_{s}(\tau_{(k)},\zeta_{(k)})dW_{s}\;,\quad \tau_{k}\wedge T \leq t < \tau_{k+1}\wedge T\;)\label{3EDSRBk-}
\enq
 for each $k$ $=$ $0,\ldots,n-1$.
\vspace{2mm}

\noindent\textbf{Step 2.} We introduce a family of processes $(\tilde{\bar Y }^k)_{0 \leq k \leq n}$ (resp. $(\tilde{\underline Y} ^k)_{0 \leq k \leq n}$). We define it recursively by 
\beqs
\tilde{\bar Y }^n_{t} & = & {\bar Y }^n_{t}(\tau_{(n)},\zeta_{(n)})\mathds{1}_{t\geq \tau_{n}} \textrm{ (resp. }\tilde{\underline Y }^n_{t} ~ = ~ {\underline Y }^n_{t}(\tau_{(n)},\zeta_{(n)})\mathds{1}_{t\geq \tau_{n}} \textrm{)}\;, \quad0\leq t\leq T\;,
\enqs
and for $k$ $=$ $0,\ldots,n-1$
\beqs
\tilde{\bar Y }^k_{t} & = & {\bar Y }^k_{t}(\tau_{(k)},\zeta_{(k)})\mathds{1}_{\tau_{k}\leq t< \tau_{k+1}} + \tilde{\bar Y }^{k+1}_{t}\mathds{1}_{t\geq \tau_{k+1}}\\
\textrm{ (resp. }\tilde{\underline Y }^k_{t}  & =  &  {\underline Y }^k_{t}(\tau_{(k)},\zeta_{(k)})\mathds{1}_{\tau_{k}\leq t< \tau_{k+1}}+ \tilde{\underline Y }^{k+1}_{t}\mathds{1}_{t\geq \tau_{k+1}} \textrm{)}\;, \quad0\leq t\leq T\;.
\enqs
These processes are c\`ad-l\`ag with jumps only at times $\tau_{l}$, $l= 1,\ldots,n$. Notice also that  $\tilde{\bar Y }^n$ (resp.  $\tilde{\underline Y }^n$, $\tilde{\bar Y }^k$, $\tilde{\underline Y }^k$) satisfies equation \reff{3EDSRBn+} (resp. \reff{3EDSRBn-}, \reff{3EDSRBk+}, \reff{3EDSRBk-}).

\vspace{2mm}
\noindent\textbf{Step 3.} We prove by a backward induction that  $\tilde {\underline Y}^{n}$ $\leq $ $\tilde{\bar Y}^{n}$ on $[\tau_{n} \wedge T ,T]$ and $\tilde{\underline Y}^{k}$ $\leq $ $\tilde {\bar Y}^{k}$ on $[\tau_{k}\wedge T, \tau_{k+1}\wedge T)$, for each $k$ $=$ $0,\ldots,n-1$. 

$\bullet$ Since $\underline \xi$ $\leq$ $\bar \xi$, $\underline F ^n$ $\leq$ $\bar  F ^n$ and $\bar F^n$ or $\underline F^n$ satisfy a comparison theorem for Brownian BSDEs, we immediately get from \reff{3EDSRBn+} and \reff{3EDSRBn-}
\beqs
\tilde{\underline  Y}^{n}_{t} & \leq & \tilde{ \bar  Y}^{n}_{t}\;, \quad \tau_{n}\wedge T \leq t\leq T\;.
\enqs

\vspace{2mm}

$\bullet$ Fix $k$ $\leq$ $n-1$ and suppose that  $\tilde{\underline Y}^{k+1}_{t}$ $\leq $ $\tilde{\bar Y}_{t}^{k+1}$ for $t$ $\in$ $[\tau_{k+1}\wedge T,\tau_{k+2}\wedge T)$. 
Denote by $^p\tilde{\bar Y}^l$ (resp. $^p\tilde{\underline Y}^l$) the predictable projection of  $\tilde{\bar Y}^l$ (resp.  $\tilde{\underline Y}^l$) for $l$ $=$ $0,\ldots,n$. Since  the random measure $\mu$ admits an intensity absolutely continuous w.r.t. the Lebesgue measure on $[0,T]$, $\tilde{\bar Y}^l$ (resp.  $\tilde{\underline Y}^l$) has inaccessible jumps (see Chapter IV of \cite{demey75}). We then have 
\beqs
^p\tilde{\bar Y}^l_{t} & = & \tilde{\bar Y}^l_{t- }\quad \textrm{(resp. }~^p\tilde{\underline Y}^l_{t} ~ = ~ \tilde{\underline Y}^l_{t- } {\rm ) }\;,\quad 0\leq t\leq T\;.
\enqs
From equations \reff{3EDSRBk+} and \reff{3EDSRBk-}, and the definition of $\tilde{\bar Y}^l$ (resp.  $\tilde{\underline Y}^l$), we have for $l$ $=$ $k$ 
\beq\nonumber
^p\tilde{\bar Y}^k_{t} & = & ^p\tilde{\bar  Y}^{k+1}_{\tau_{k+1}} \mathds{1}_{\tau_{k+1} \leq T} + \bar \xi \mathds{1}_{ \tau_{k+1} > T}
+\int_{t}^{\tau_{k+1}\wedge T} \bar F^k\Big(s,^p\tilde{\bar Y}^k_{s},\bar Z^k_{s}(\tau_{(k)},\zeta_{(k)})\Big)ds\qquad\\
& &  -\int_{t}^{\tau_{k+1}\wedge T}\bar Z^k_{s}(\tau_{(k)},\zeta_{(k)})dW_{s}\;,\quad \tau_{k} \wedge T \leq t < \tau_{k+1}\wedge T\;.\\\nonumber
\\\nonumber
({\rm resp. }\quad ^p\tilde{\underline  Y}^k_{t} & = & ^p\tilde{\underline Y}^{k+1}_{\tau_{k+1}} \mathds{1}_{\tau_{k+1} \leq T} + \underline \xi \mathds{1}_{ \tau_{k+1} > T}+\int_{t}^{\tau_{k+1}\wedge T} \underline F^k\Big(s,^p\tilde{\underline Y}^k_{s},\underline Z^k_{s}(\tau_{(k)},\zeta_{(k)})\Big)ds\qquad\nonumber\\
& &  -\int_{t}^{\tau_{k+1} \wedge T}\underline Z^k_{s}(\tau_{(k)},\zeta_{(k)})dW_{s}\;,\quad \tau_{k}\wedge T \leq t < \tau_{k+1}\wedge T\;) 
\enq
Since $\tilde{\bar  Y}^{k+1}_{\tau_{k+1}}$ $\geq$ $\tilde{\underline Y}^{k+1}_{\tau_{k+1}}$,  we get $^p\tilde{\bar  Y}^{k+1}_{\tau_{k+1}}$ $\geq$ $^p\tilde{\underline Y}^{k+1}_{\tau_{k+1}}$. This together with conditions on $\bar \xi$, $\underline \xi$, $\bar F ^k$ and $\underline F ^k$ give the result. 
\vspace{2mm}

\noindent\textbf{Step 4.} Since $\tilde{\bar Y}^k$ (resp. $\tilde{\underline Y}^k$) coincides with $\bar Y$ (resp. $\underline Y$) on $[\tau_{k}\wedge T,\tau_{k+1}\wedge T)$, we get the result.\ep
\begin{Remark}
{\rm   It is possible to obtain Theorem 4.1 under weaker assumptions than (HC). For instance, it is sufficient to assume that $W$ is a $\G$-semimartingale fo the form
\beqs
W & = & M+\int_0^.a_sds\;,
\enqs
 with $M$ a $\G$-local martingale and $a$ a $\G$-adapted process satisfying 
\beq\label{cond-girsanov}
\E\Big[\exp\Big(-\int_0^Ta_sdM_s-{1\over 2}\int_0^T|a_s|^2ds\Big)\Big] & = & 1\;.
\enq
Indeed, we first notice that $(M_t)_{t\in[0,T]}$ is a $\G$-Brownian motion since it is a continuous $\G$-martingale with $\langle M, M\rangle_t =t$ for $t\geq 0$.
Then, from \reff{cond-girsanov}  we can apply Girsanov Theorem and get that $(W_t)_{t\in[0,T]}$  is a $(\Q,\G)$-Brownian motion where  $\Q$ is the probability measure equivalent to $\P$ defined by
\beqs
\left.{d \Q \over d\P}\right|_{\Gc_T} & = & \exp\Big(-\int_0^Ta_sdM_s-{1\over 2}\int_0^T|a_s|^2ds\Big)\;.
\enqs
Therefore we can prove Theorem \ref{3theoreme comparaison lipschitz} under $\Q$. 
Since $\Q$ is equivalent to  $\P$ the conclusion remains true under $\P$. 
%
}
\end{Remark}

\subsection{Uniqueness via comparison}
In this form, the previous theorem is not usable since the condition on the generators of the Brownian BSDEs is implicit: it involves the solution of the previous Brownian BSDEs at each step. 
We give, throughout the sequel, an explicit example for which Theorem \ref{3theoreme comparaison lipschitz} provides uniqueness. This example is based on a comparison theorem for quadratic BSDEs given by Briand and Hu \cite{bh08}.  We first introduce the following assumptions. 
\begin{enumerate}[(HUQ1)]
\item The function $f(t,y,.,u)$ is concave for all $(t,y,u)\in[0,T]\times\R\times Bor(E,\R)$.
\item There exists a constant $L$ s.t. 
\beqs
|f(t,y,z,(u(e)-y)_{e\in E})-f(t,y',z,(u(e)-y')_{e\in E})| & \leq & L|y-y'| 
\enqs
for all $(t,y,y',z ,u )\in[0,T]\times[\R]^2\times\R^d\times Bor(E,\R)$. 
\item There exists a constant $C>0$ such that 
\beqs
|f(t,y,z,u)| & \leq & C \Big(1+|y|+|z|^2+\int_E|u(e)|\lambda_t(e)de\Big)
\enqs 
for all $(t,y,z,u)\in[0,T]\times\R\times\R^d\times Bor(E,\R)$. 
\item $f(t,.,u)=f(t,.,0)$ for all $u\in Bor(E,\R)$ and all $t\in(\tau_n\wedge T,T]$.
\end{enumerate} 



\begin{Theorem}\label{TUQ} Under (HD), (HBI), (HC), (HUQ1), (HUQ2), (HUQ3) and (HUQ4), 
BSDE \reff{3BSDE jump} admits at most one solution.
\end{Theorem}
\noindent\textbf{Proof.} Let $(Y,Z,U)$ and $(Y',Z',U')$ be two solutions of \reff{3BSDE jump} in  $\Sc_\G^\infty[0,T]\times L^2_{\G}[0,T]\times L^2(\mu)$. Define the process $\tilde U$ (resp. $\tilde U'$) by
\beqs
\tilde U_t(e)~~ (resp. \quad\tilde U'_t(e)) & = & U_t(e)\mathds{1}_{t \leq \tau_n} ~~(resp. \quad U'_t(e)\mathds{1}_{t \leq \tau_n}) \;,\quad (t,e)\in [0,T]\times E\;.
\enqs
  Then, $U=\tilde U$ and $ U'= \tilde U'$ in $L^2(\mu)$. Therefore, from (HUQ4), $(Y,Z,\tilde U)$ and $(Y',Z',\tilde U')$ are also solutions to \reff{3BSDE jump} in  $\Sc_\G^\infty[0,T]\times L^2_{\G}[0,T]\times L^2(\mu)$.
  
  We now prove by a backward induction on $k=n,n-1,\ldots,1,0$ that 
  \beqs
  Y_t\mathds{1}_{\tau_k\leq t} & = & Y'_t\mathds{1}_{\tau_k\leq t}\;, \quad t\in[0,T]\;,
  \enqs
    $\bullet$ Suppose that $k=n$.  Then, $(Y_t\mathds{1}_{\tau_n\leq t},Z_t\mathds{1}_{\tau_n< t},(\tilde U_t+Y_{t^-})\mathds{1}_{\tau_{n-1}< t\leq\tau_n})$ and $(Y'_t\mathds{1}_{\tau_n\leq t},Z'_t\mathds{1}_{\tau_n< t},(\tilde U'_t+Y_{t-})\mathds{1}_{\tau_{n-1}< t\leq\tau_n})$ are solution to
  \beqs
  Y_t & = & \xi \mathds{1}_{\tau_n\leq T}+ \int_t^T\mathds{1}_{\tau_n< s}f(s,Y_s,Z_s,0)ds-\int_t^TZ_s dW_s-\int_t^T\int_E U_s(e)\mu(de,ds)\;,\quad t\in[0,T]\;.
  \enqs
 Using Remark \ref{rem W G-MG} and Theorem 5 in \cite{bh08}, we obtain that the generator $\mathds{1}_{\tau_n < .}f$ satisfies a comparison theorem in the sense of Definition \ref{defcompB}.  We can then apply Theorem \ref{3theoreme comparaison lipschitz} with 
\beqs
\underline F(t,y,z,u)  & = & \bar F(t,y,z,u) ~=~  \mathds{1}_{\tau_n < t}f(t,y,z,0)\;,\quad (t,y,z,u)\in[0,T]\times\R\times\R^d\times Bor(E,\R)\;,
\enqs
 and we get that $Y_.\mathds{1}_{ \tau_n \leq .}= Y'_.\mathds{1}_{ \tau_n \leq .}$.   
  \vspace{1mm}

 \ni  $\bullet$ Suppose that  $Y_.\mathds{1}_{ \tau_{k+1} \leq . }= Y'_.\mathds{1}_{\tau_{k+1} \leq . }$.  We can then choose $Y^{j}$
  and $Y'^{j}$ appearing in the decomposition of the processes $Y$ and $Y'$ given by Lemma \ref{3lemme decomposition} (ii) such that
  \beqs
  Y^{j}_s(\theta_{(j)},e_{(j)}) & = & Y'^{j}_s(\theta_{(j)},e_{(j)})\;,
  \enqs
  for all $(\theta,e)\in \Delta_n\times E^n$ and $j= k+1, \ldots, n$. 
 Therefore, we get that $(Y_t\mathds{1}_{\tau_k\leq t},Z_t\mathds{1}_{\tau_k< t},(\tilde U_t+Y_{t^-}\mathds{1}_{ t\leq\tau_k})\mathds{1}_{\tau_{k-1}< t})$ and $(Y'_t\mathds{1}_{\tau_k\leq t},Z'_t\mathds{1}_{\tau_k< t},(\tilde U'_t+Y_{t^-}\mathds{1}_{ t\leq\tau_k})\mathds{1}_{\tau_{k-1}< t})$ are solution to
  \beqs
  Y_t & = & \xi \mathds{1}_{\tau_k\leq T} + \int_t^TF(s,Y_s,Z_s)ds-\int_t^TZ_s dW_s-\int_t^T\int_E U_s(e)\mu(de,ds)\;,\quad t\in[0,T]\;,
\enqs
where the generator $F$ is defined by
\beqs
F(t,y,z) & = & \sum_{j=k}^{n-1}\mathds{1}_{\tau_k<t\leq \tau_{k+1}}F^k(t,y,z)+\mathds{1}_{\tau_n<t}F^n(t,y,z)
\;,
\enqs
where 
\beqs
F^k(t,y,z) & =& f\Big(t,y,z,Y^{k+1}_s(\tau_{(k)},s,\zeta_{(k)},.)-y,\tau_{(k)},\zeta_{(k)}\Big)
\enqs
\beqs
F^n(t,y,z) & = &  f(t,y,z,0)
\enqs
 for all $(t,y,z)\in[0,T]\times\R\times\R^d$.
Using Remark \ref{rem W G-MG} and Theorem 5 in \cite{bh08}, we obtain that the generator $F$ satisfies a comparison theorem in the sense of Definition \ref{defcompB}. 
We can then apply Theorem \ref{3theoreme comparaison lipschitz} and we get that $Y_.\mathds{1}_{\tau_k \leq .} =Y'\mathds{1}_{\tau_k \leq .} $. 

\vspace{1mm}

\ni $\bullet$ Finally  the result holds true for all $k=0,\ldots,n$ which gives  $Y=Y'$.

\vspace{1mm}

\ni $\bullet$ We now prove that $Z=Z'$ and $ U=U'$. 
Identifying the finite variation part and the unbounded variation part of $Y$ we get $Z=Z'$. Then, identifying the pure jump part of $Y$ we get $\tilde U =\tilde U'$. Since $\tilde U =U$ (resp. $\tilde U' =U'$ ) in $L^2(\mu)$, we finally get $(Y,Z,U)=(Y',Z',U')$. 
\ep

\section{Exponential utility maximization in a jump market model}\label{UtilMaxjump}
\setcounter{equation}{0} \setcounter{Assumption}{0}
\setcounter{Theorem}{0} \setcounter{Proposition}{0}
\setcounter{Corollary}{0} \setcounter{Lemma}{0}
\setcounter{Definition}{0} \setcounter{Remark}{0}

We consider a financial market model with a riskless bond assumed for simplicity equal to one, and a risky asset subjects to some counterparty risks. We suppose that the Brownian motion $W$ is one dimensional ($d=1$). The dynamic of the risky asset is affected by other firms, the counterparties, which may default at some random times, inducing consequently some jumps in the asset price. However, this asset still exists and can be traded after the default of the counterparties. We keep the notation of previous sections.\\
\indent Throughout the sequel, we suppose that (HD), (HBI) and (HC) are satisfied. 
We consider that the price process $S$ evolves according to the equation 
\beqs
S_t=S_0+\int_0^tS_{u^-}\Big(b_udu+\sigma_udW_u+\int_{E}\beta_u(e)\mu(de,du)\Big)\;,\quad 0 \leq t \leq T\;.
\enqs
All processes $b$, $\sigma$ and $\beta$ are assumed to be $\G$-predictable. We introduce the following assumptions on the coefficients appearing in the dynamic of $S$:

\vspace{2mm}
\begin{enumerate}[(HS1)]
\item The processes $b$, $\sigma$ and $\beta$ are uniformly bounded: there exists a constant $C$ s.t.
\beqs
|b_t|+|\sigma_t|+|\beta_t(e)| & \leq & C\;,~~0 \leq t \leq T\;, ~e \in E\;,~~ \P-a.s.
\enqs

\item There exists a positive constant $c_\sigma$ such that
\beqs
\sigma_t \geq c_\sigma \;, \quad 0 \leq t \leq T\;, \quad \P-a.s. 
\enqs

\item The process $\beta$ satisfies:
\beqs
\beta_t(e) >-1\;, \quad 0 \leq t \leq T\;,\quad e\in E\;, \quad \P-a.s. 
\enqs
\item The process $\vartheta$ defined by $
\vartheta_t  =  \frac{b_t}{\sigma_t}$, $t\in [0,T] $,
is uniformly bounded: there exists a constant $C$ such that
\beqs
|\vartheta_t| & \leq & C\;,~~ 0 \leq t \leq T\;, ~~Ê \P-a.s.
\enqs
\end{enumerate}

\vspace{2mm}

\ni We notice that (HS1) allows the process $S$ to be well defined and (HS3) ensures it to be positive. 


\vspace{2mm}

A self-financing trading strategy is determined by its initial capital $x \in \R$ and the amount of money $\pi_t$ invested in the stock, at time $t \in [0, T]$. 
The wealth at time $t$ associated with a strategy $(x, \pi)$ is
\beqs
X^{x,\pi}_t & = &x+\int_0^t\pi_s b_sds+\int_0^t\pi_s\sigma_sdW_s+\int_0^t\int_E\pi_s\beta_s(e)\mu(de,ds)\;,~~0 \leq t \leq T\;.
\enqs
We consider a contingent claim, that is a random payoff at time $T$ described by a $\Gc_T$-measurable random variable $B$. We suppose that $B$ is bounded and satisfies
\beqs
B & = & \sum_{k=0}^{n}B^k(\tau_{(k)}, \zeta_{(k)})\mathds{1}_{\tau_k\leq T< \tau_{k+1}} \;,
\enqs
where $B^0$ is $\Fc_{T}$-measurable and $B^k$ is $\Fc_{T}\otimes\Bc(\Delta_{k})\otimes\Bc(E^{k})$-measurable for each $k=1,\ldots,n$. Then, we define
\beq\label{3V}
V(x) & = & \sup\limits_{\pi \in \Ac} \mathbb{E}\big[- \exp( - \alpha (X_T^{x,\pi}-B))\big] \;,
\enq
the maximal expected utility that we can achieve by starting at time $0$ with the initial capital $x$, using some admissible strategy $\pi \in \Ac$ (which is defined throughout the sequel) on $[0,T]$ and paying $B$ at time $T$. $\alpha$ is a given positive constant which can be seen as a coefficient of absolute risk aversion.

Finally, we introduce a compact subset $\Cc$ of $\R$ with $0 \in \Cc$, which represents an eventual  constraint imposed to the trading strategies, that is, $\pi_t(\omega)$ $\in$ $\Cc$.  
We then define the space $\Ac$ of admissible strategies. 
\begin{Definition}\label{3ensemble admissible}{\rm
The set $\Ac$ of admissible strategies consists of all $\R$-valued $\Pc(\G)$-measurable processes $\pi$  $=$ $(\pi_t)_{0\leq t\leq T}$ which satisfy $\E\int_0^T|\pi_t\sigma_t|^2dt + \E \int_0^T \int_E |\pi_t \beta_t(e)| \lambda_t(e) de dt$ $<$ $\infty$, and $\pi_t$ $\in$ $\Cc$, $dt\otimes d\P-a.e.$, as well as  
the uniform integrability of the family
\beqs
\Big\{ \exp \big(-\alpha X^{x,\pi}_\tau\big)~:~\tau \text{ stopping time valued in }[0,T] \Big\}\;.
\enqs 
 }
\end{Definition}


We first notice that the compactness of $\Cc$ implies the integrability conditions imposed to the admissible strategies. 

\begin{Lemma}\label{lem ttes adm}
Any $\Pc(\G)$-measurable process $\pi$ valued in $\Cc$ satisfies $\pi\in\Ac$.
\end{Lemma}
\ni The proof is exactly the same as in \cite{mor09}. We therefore omit it.

\vspace{2mm}

In order to characterize the value function $V(x)$ and an optimal strategy, we construct, as in \cite{huimkmul05} and \cite{mor09}, a family of stochastic processes $(R^{(\pi)})_{\pi \in \Ac}$ with the following properties:
\begin{enumerate}[(i)]
\item $R^{(\pi)}_T=-\exp(-\alpha(X^{x,\pi}_T- B))$ for all $\pi\in\Ac$, 
\item $R^{(\pi)}_0=R_0$ is constant for all $\pi\in \Ac$,
\item $R^{(\pi)}$ is a supermartingale for all $\pi\in \Ac$ and there exists $\hat{\pi}\in \Ac$ such that  $R^{(\hat{\pi})}$ is a martingale.
\end{enumerate}
Given processes owning these properties we can compare the expected utilities of the strategies $\pi\in \Ac$ and $\hat{\pi}\in \Ac$ by
\beqs
\E\big[-\exp \big(-\alpha(X^{x,\pi}_T- B)\big)\big]&\leq & R_0(x) ~ = ~ \E\big[-\exp\big(-\alpha(X^{x,\hat{\pi}}_T- B)\big)\big] ~ = ~ V(x) \;,
\enqs
whence $\hat{\pi}$ is the desired optimal strategy. To construct this family, we set 
\beqs
R^{(\pi)}_t & = & -\exp\big(-\alpha(X^{x,\pi}_t-Y_t)\big)\;,~~0\leq t\leq T\;,~\pi\in\Ac \;,
\enqs
where $(Y, Z, U)$ is a solution of the BSDE
\beq\label{3equation f}
Y_t & = & B + \int_t^Tf(s,Z_s,U_s)ds-\int_t^TZ_sdW_s-\int_t^T\int_{E}U_s(e)\mu(de,ds)\;,\quad 0\leq t\leq T\;. \qquad 
\enq
We have to choose a function $f$ for which $R^{(\pi)}$ is a supermartingale for all $\pi\in \Ac$, and there exists a $\hat{\pi}\in \Ac$ such that $R^{(\hat{\pi})}$ is a martingale. We assume that there exists a triple $(Y, Z, U)$ solving a BSDE with jumps of the form (\ref{3equation f}), with terminal condition $B$ and with a driver $f$ to be determined. We first apply It\^{o}'s formula to $R^{(\pi)}$ for any strategy $\pi$:
\beqs
dR^{(\pi)}_t&=&R^{(\pi)}_{t^-}\Big[\Big(-\alpha \big(f(t,Z_t,U_t)+\pi_t b_t \big)+\frac{\alpha^2}{2}(\pi_t\sigma_t-Z_t)^2\Big)dt-\alpha(\pi_t\sigma_t-Z_t)dW_t\\
&&+\int_E \big(\exp\big(-\alpha\big(\pi_t\beta_t(e)-U_t(e)\big)\big)-1 \big)\mu(de,dt)\Big].
\enqs
Thus, the process $R^{(\pi)}$ satisfies the following SDE: 
\beqs
dR^{(\pi)}_t=R^{(\pi)}_{t^-}dM^{(\pi)}_t+R^{(\pi)}_tdA^{(\pi)}_t\;, \quad 0< t\leq T\;,
\enqs
 with  $M^{(\pi)}$ a local martingale and $A^{(\pi)}$ a finite variation continuous process given by
\[\left\{
\begin{aligned}
dM^{(\pi)}_t  = & ~-\alpha(\pi_t\sigma_t-Z_t)dW_t+\int_E \big(\exp\big(-\alpha \big(\pi_t\beta_t(e)-U_t(e)\big)\big)-1 \big)\tilde{\mu}(de,dt)\,,\\
dA^{(\pi)}_t  = &~\Big(-\alpha \big(f(t,Z_t,U_t)+\pi_t b_t\big)+\frac{\alpha^2}{2}(\pi_t\sigma_t-Z_t)^2\\
 &~+\int_E \big(\exp\big(-\alpha \big(\pi_t\beta_t(e)-U_t(e)\big)\big)-1 \big)\lambda_t(e) de\Big)dt \;.
\end{aligned}
\right.\]
It follows that $R^{(\pi)}$ has the multiplicative form
\beqs
R^{(\pi)}_t=R^{(\pi)}_0\mathfrak{E}(M^{(\pi)})_t\exp\big(A^{(\pi)}_t\big) \;,
\enqs
where $\mathfrak{E}(M^{(\pi)})$ denotes the Doleans-Dade exponential of the local martingale $M^{(\pi)}$. Since $\exp(-\alpha(\pi_t\beta_t(e)-U_t(e)))-1> -1,~\P-a.s.$, the Doleans-Dade exponential of the discontinuous part of  $M^{(\pi)}$ is a positive local martingale and hence, a supermartingale. The supermartingale condition in (iii) holds true, provided, for all $\pi\in\Ac$, the process $\exp(A^{(\pi)})$ is nondecreasing, this entails
\beqs
-\alpha \big(f(t,Z_t,U_t)+\pi_t b_t\big)+\frac{\alpha^2}{2}(\pi_t\sigma_t-Z_t)^2+\int_E \big(\exp\big(-\alpha\big(\pi_t\beta_t(e)-U_t(e)\big)\big)-1\big)\lambda_t(e) de\geq 0\;.
\enqs
This condition holds true, if we define $f$ as follows
\beqs
f(t,z,u) & = & \inf_{\pi\in \Cc}\Big\{\frac{\alpha}{2}\Big|\pi\sigma_t-\Big(z+\frac{\vartheta_t}{\alpha}\Big)\Big|^2 + \int_E\frac{\exp(\alpha(u(e)-\pi\beta_t(e)))-1}{\alpha}\lambda_t(e)de\Big\} \\
&&- \vartheta_tz - \frac{|\vartheta_t|^2}{2\alpha}\;,
\enqs
recall that $\vartheta_{t}$ $=$ $b_{t} / \sigma_{t}$ for $t$ $\in$ $[0,T]$.
\begin{Theorem}
Under (HD), (HBI), (HC), (HS1), (HS2), (HS3) and (HS4), the value function of the optimization problem (\ref{3V}) is given by 
\beq\label{3fonction valeur}
V(x) & = & -\exp(-\alpha(x-Y_0))\,,
\enq
where $Y_0$ is defined as the initial value of the unique solution $(Y,Z,U)\in \Sc^\infty_{\G}[0,T]\times L^2_{\G}[0,T]\times L^2(\mu)$ of the BSDE
\beq\label{3equation exponentielle}
Y_t & = & B + \int_t^Tf(s,Z_s,U_s)ds-\int_t^TZ_sdW_s-\int_t^T\int_{E}U_s(e)\mu(de,ds)\;, \quad 0\leq t\leq T\;,\qquad 
\enq
with  
\beqs
f(t,z,u) & = & \inf_{\pi\in \Cc}\Big\{\frac{\alpha}{2}\Big|\pi\sigma_t-\Big(z+\frac{\vartheta_t}{\alpha}\Big)\Big|^2+\int_E\frac{\exp(\alpha(u(e)-\pi\beta_t(e)))-1}{\alpha}\lambda_t(e) de\Big\} \\
&&- \vartheta_tz - \frac{|\vartheta_t|^2}{2\alpha}\;,
\enqs
for all $(t,z,u)\in[0,T]\times\R\times Bor(E,\R)$.
There exists an optimal trading strategy $\hat{\pi} \in \Ac$ which satisfies
\beq\label{pi optimal}
\hat{\pi}_{t} &\in & \emph{arg}\min_{\pi\in \Cc}\Big\{\frac{\alpha}{2}\Big|\pi\sigma_t-\Big(z+\frac{\vartheta_t}{\alpha}\Big)\Big|^2+\int_E\frac{\exp(\alpha(u(e)-\pi\beta_t(e)))-1}{\alpha}\lambda_t(e) de\Big\} \;,
\enq
for all $t$ $\in$ $[0,T]$.
\end{Theorem}

\noindent\textbf{Proof.}
\textbf{Step 1.} We first prove the existence of a solution to BSDE \reff{3equation exponentielle}.  We first check the measurability of the generator $f$.
Notice that we have $f(.,.,.,.) =\inf_{\pi\in\Cc }F(\pi,.,.,.,.)$ where $F$ is defined by
\beqs
F(\pi,t,y,z,u) & = & \frac{\alpha}{2}\Big|\pi\sigma_t-\Big(z+\frac{\vartheta_t}{\alpha}\Big)\Big|^2+\int_E\frac{\exp(\alpha(u(e)-\pi\beta_t(e)))-1}{\alpha}\lambda_t(e) de
\enqs
for all $(\omega,t,\pi,y,z,u)\in\Omega\times[0,T]\times\Cc\times\R\times\R \times Bor(E,\R)$. From Fatou's Lemma we have that $u\mapsto\int_Eu(e)de$ is l.s.c. and hence measurable on $Bor(E,\R_+):=\{u\in Bor(E,\R)~:~u(e)\geq0\;,~\forall e\in E\}$. Therefore $F(\pi,.,.,.,.)$ is $\Pc(\G)\otimes\Bc(\R)\otimes\Bc(\R) \otimes\Bc(Bor(E,\R))$-measurable for all $\pi\in\Cc$. 
Since $F(.,t,y,z,u)$ is continuous for all $(t,y,z,u)$ we have
$f(.,.,.,.) =\inf_{\pi\in\Cc\cap\Q }F(\pi,.,.,.,.)$, and $f$ is $\Pc(\G)\otimes\Bc(\R)\otimes\Bc(\R)\otimes\Bc(Bor(E,\R))$-measurable.

\vspace{1mm}

We now  apply Theorem \ref{3existence solution}.  
Let $\sigma^k$, $\vartheta^k$ and $\beta^k$, $k$ $=$ $0,\ldots,n$, be the respective terms appearing in the decomposition  of $\sigma$, $\vartheta$ and $\beta$  given by Lemma \ref{3lemme decomposition}. Using  (HS1) and (HS4), we can assume w.l.o.g. that these terms are uniformly bounded. 
Then, in the decomposition of the generator $f$, we can choose the functions $f^k$, $k=0,\ldots,n$, as
\beqs
f^n(t,z,u,\theta,e) & = & \inf_{\pi\in \Cc}\Big\{\frac{\alpha}{2}\Big|\pi\sigma^n_t(\theta,e)-\Big(z+\frac{\vartheta^n_t(\theta,e)}{\alpha}\Big)\Big|^2\Big\} - \vartheta^n_t(\theta,e)z - \frac{|\vartheta^n_t(\theta,e)|^2}{2\alpha}\;,
\enqs
and
\beqs
f^k(t,z,u,\theta_{(k)},e_{(k)}) & = & \inf_{\pi\in \Cc}\Big\{\frac{\alpha}{2}\Big|\pi\sigma^k_t(\theta_{(k)},e_{(k)})-\Big(z+\frac{\vartheta^k_t(\theta_{(k)},e_{(k)})}{\alpha}\Big)\Big|^2\\
& & +\int_E\frac{\exp(\alpha(u(e')-\pi\beta^k_t(\theta_{(k)},e_{(k)},e')))-1}{\alpha}\lambda^{k+1}_t(e',\theta_{(k)},e_{(k)}) de'\Big\} \\
&&- \vartheta^k_t(\theta_{(k)},e_{(k)})z - \frac{|\vartheta^k_t(\theta_{(k)},e_{(k)})|^2}{2\alpha}\;,
\enqs
 for $k=0,\dots,n-1$ and $(\theta,e)\in \Delta_n\times E^n$.
 
 Notice also that since $B$ is bounded, we can choose $B^k$, $k=0,\ldots,n$,  uniformly bounded. 
We now prove by backward induction on $k$ that the BSDEs (we shall omit the dependence on $(\theta,e)$)
\beq \label{3BSDE n}
Y_{t}^n & = & B^n +\int_{t}^Tf^n(s,Z^n_{s},0)ds -\int_{t}^TZ^n_{s}dW_{s}\;, \quad \theta_n\wedge T\leq t\leq T\;,~(k=n)
\enq
and 
\beq\nonumber
Y_{t}^k & = & B^k +\int_{t}^Tf^k(s,Z^k_{s},Y^{k+1}_{s}(s,.)-Y^k_{s})ds\\
 & &  -\int_{t}^TZ^k_{s}dW_{s}\;, \quad \theta_k\wedge T\leq t\leq T\;,~(k=0,\ldots,n-1) \label{3BSDE k}
\enq
admit a solution  $(Y^k,Z^k)$ in $\Sc^\infty_\F[\theta_k\wedge T,T] \times L^2_{\F}[\theta_k\wedge T,T]$ such that $Y^k$ (resp. $Z^k$) is $\Pc\Mc(\F)\otimes\Bc(\Delta_k)\otimes\Bc(E^k)$ (resp. $\Pc(\F)\otimes\Bc(\Delta_k)\otimes\Bc(E^k)$)-measurable 
with 
\beqs
\sup_{(\theta,e)\in\Delta_{n}\times E^n}\|Y^k(\theta_{(k)},e_{(k)})\|_{\Sc^\infty[\theta_{k}\wedge T,T]} +\|Z^k(\theta_{(k)},e_{(k)})\|_{L^2[\theta_{k}\wedge T,T]}
& < & \infty\;, 
\enqs
for all $k$ $=$ $0,\ldots,n$.

\ni$\bullet$ Since $0$ $\in$ $\Cc$, we have 
\beqs
- \vartheta^n_t z - \frac{|\vartheta^n_t|^2}{2\alpha} \leq  & f^n(t,z,0)  \leq & \frac{\alpha}{2}|z|^2 \;.
\enqs
Therefore, we can apply Theorem 2.3 of \cite{kob00}, and we get that for any $(\theta,e)\in\Delta_{n}\times E^n$, there exists a solution $\big(Y^n(\theta,e),Z^n(\theta,e)\big)$ to BSDE (\ref{3BSDE n}) in $\Sc^\infty_\F[\theta_n\wedge T,T]\times L^2_\F[\theta_n\wedge T,T]$. Moreover, this solution is constructed as a limit of Lipschitz BSDEs (see \cite{kob00}). Using Proposition \ref{mesLipsol}, we get that $Y^n$ (resp. $Z^n$) is $\Pc\Mc(\F)\otimes\Bc(\Delta_n)\otimes\Bc(E^n)$ (resp. $\Pc(\F)\otimes\Bc(\Delta_n)\otimes\Bc(E^n)$)-measurable. 

Then, using Proposition 2.1 of \cite{kob00}, we get the existence of a constant $K$ such that 
\beqs
\sup_{(\theta,e)\in\Delta_{n}\times E^n}{\|Y^n(\theta,e)\|}_{\Sc^\infty[\theta_n\wedge T,T]} +{\|Z^n(\theta,e)\|}_{L^2[\theta_n\wedge T,T]}
& \leq & K\;.
\enqs
$\bullet$ Suppose that BSDE \reff{3BSDE k} admits a solution at rank $k+1$ ( $k \leq n-1$) with
\beq\nonumber
\sup_{(\theta,e)\in\Delta_{n}\times E^n}\Big\{{\|Y^{k+1}(\theta_{(k+1)},e_{(k+1)})\|}_{\Sc^\infty[\theta_{k+1}\wedge T,T]} \quad & & \\
 \quad  +{\|Z^{k+1}(\theta_{(k+1)},e_{(k+1)})\|}_{L^2[\theta_{k+1}\wedge T,T]}\Big\}\label{hypreck+1}
& < & \infty\;. 
\enq
We denote $g^k$ the function defined by
\beqs
g^k(t,y,z,\theta_{(k)},e_{(k)}) & = & f^k(t,z,Y^{k+1}_{t}(\theta_{(k)},t,e_{(k)},.)-y,\theta_{(k)},e_{(k)})\;,
\enqs
for all $(t,y,z)\in[0,T]\times\R\times\R$ and $(\theta,e)\in\Delta_n\times E^n$. Since $g^k$ has an exponential growth in the variable $y$ in the neighborhood of $-\infty$, we can not directly apply our previous results. We then prove via a comparison theorem that there exists a solution by introducing another BSDE  which admits a solution and whose generator coincides with $g$ in the domain where the solution lives.
  
Let $(\underline Y^k(\theta_{(k)},e_{(k)}), \underline Z^k(\theta_{(k)},e_{(k)}))$ be the solution in $\Sc^\infty_\F[\theta_{k}\wedge T,T]\times L^2_\F[\theta_{k}\wedge T,T]$ to the linear BSDE
\beqs
\underline Y^k_t (\theta_{(k)},e_{(k)})& = & B^k(\theta_{(k)},e_{(k)}) +\int_t^T\underline g^k(s,\underline Y^k_s,\underline Z^k_s)(\theta_{(k)},e_{(k)})ds\\
 & & -\int_t^T\underline Z^k_s(\theta_{(k)},e_{(k)})dW_s\;,\qquad \theta_k\wedge T\leq t\leq T\;, 
\enqs 
where 
\beqs
\underline g^k(t,y,z,\theta_{(k)},e_{(k)}) & = & -\vartheta_t^k(\theta_{(k)},e_{(k)}) z-\frac{\vartheta^k_t(\theta_{(k)},e_{(k)})}{2\alpha}\;, 
\enqs
for all $(t,y,z)\in[0,T]\times\R\times \R$.
Since $B^k$ and $\vartheta^k$ are uniformly bounded, we have
\beq\label{yunderline}
\sup_{(\theta_{(k)},e_{(k)}) \in \Delta_k \times E^k}{\|\underline Y^{k}(\theta_{(k)},e_{(k)})\|}_{\Sc^\infty[\theta_k\wedge T,T]} & < & \infty\;.
\enq
Then, define the generator $\tilde g^k$ by 
\beqs
\tilde g^k(t,y,z,\theta_{(k)},e_{(k)}) & = & g^k(t,y\vee \underline Y^k_t(\theta_{(k)},e_{(k)}),z,\theta_{(k)},e_{(k)})\;, 
\enqs
 for all $(t,y,z)\in[0,T]\times\R\times \R$ and $(\theta,e)\in\Delta_n\times E^n$. 

Moreover, since $0$ $\in$ $\Cc$, 
we get from \reff{hypreck+1} and \reff{yunderline} the existence of a positive constant $C$ such that 
\beqs
|\tilde g^k(t,y,z,\theta_{(k)},e_{(k)})| & \leq & C(1+|z|^2)\;,
\enqs
for all $(t,y,z)\in[0,T]\times\R\times\R$ and $(\theta,e)\in\Delta_n\times E^n$. 
We can then apply Theorem 2.3 of \cite{kob00}, and we obtain that the BSDE
\beqs
\tilde Y^k_t(\theta_{(k)},e_{(k)}) & = & B^k(\theta_{(k)},e_{(k)}) +\int_t^T\tilde g ^k(s ,\tilde Y_s^k,\tilde Z_s^k)(\theta_{(k)},e_{(k)})ds\\
 & & -\int_t^T\tilde Z_s^k(\theta_{(k)},e_{(k)}) dW_s\;,\quad \theta_k\wedge T\leq t\leq T \;,
\enqs
 admits a solution $(\tilde Y^k(\theta_{(k)},e_{(k)}),\tilde Z^k(\theta_{(k)},e_{(k)}))\in \Sc^\infty_\F[\theta_k\wedge T,T] \times L^2_{\F}[\theta_k\wedge T,T]$. Using Proposition 2.1 of \cite{kob00}, we get  
 \beqs
 \sup_{(\theta_{(k)},e_{(k)}) \in \Delta_k \times E^k}\|\tilde Y^{k}(\theta_{(k)},e_{(k)})\|_{\Sc^\infty[\theta_k\wedge T,T]} & < & \infty \;.
 \enqs
Then, since $\tilde g^k$ $\geq$ $\underline g^k$ and since $\underline g^k$ is Lipschitz continuous, we get from the comparison theorem for BSDEs that 
$\tilde Y^k$ $\geq$ $\underline Y^k$. Hence, $(\tilde Y^k,\tilde Z^k)$ is solution to BSDE \reff{3BSDE k}.
Notice then that we can choose $\tilde Y^k$ (resp. $\tilde Z^k$) as a $\Pc\Mc(\F)\otimes\Bc(\Delta_k)\otimes\Bc(E^k)$ (resp. $\Pc(\F)\otimes\Bc(\Delta_k)\otimes\Bc(E^k)$)-measurable process. Indeed, these processes are solutions to quadratic BSDEs and hence can be written as the limit of solutions to Lipschitz BSDEs (see \cite{kob00}). Using Proposition \ref{mesLipsol} with $\Xc=\Delta_k \times E^k$ and $d\rho(\theta,e)=\gamma_0(\theta,e)d\theta de$ we get that the solutions to Lipschitz BSDEs are  $\Pc(\F)\otimes\Bc(\Delta_k)\otimes\Bc(E^k)-$measurable and hence $\tilde Y^k$ (resp. $\tilde Z^k)$ is $\Pc\Mc(\F)\otimes\Bc(\Delta_k)\otimes\Bc(E^k)$ (resp. $\Pc(\F)\otimes\Bc(\Delta_k)\otimes\Bc(E^k)$)-measurable.
\\

\noindent\textbf{Step 2.} We now prove the uniqueness of a solution to BSDE \reff{3equation exponentielle}. Let $(Y^1,Z^1,U^1)$ and $(Y^2,Z^2,U^2)$ be two solutions of BSDE \reff{3equation exponentielle} in 
$\Sc_{\G}^\infty[0,T]\times L_{\G}^2[0,T]\times L^2(\mu)$. 

Applying an exponential change of variable, we obtain that $(\tilde Y^i, \tilde  Z^i, \tilde U^i)$ defined for $i=1,2$  by
\beqs
\tilde Y^i_t & = &  \exp(\alpha Y^i_t)\;,\\
\tilde Z^i_t & = & \alpha  \tilde Y^i_tZ^i_t\;,\\
\tilde U^i_t(e) & = &  \tilde Y^i_{t^-}\big( \exp(
\alpha U^i_t(e))-1 \big)\;,
\enqs
for all $t\in[0,T]$,  
are solution in $\Sc_{\G}^\infty[0,T]\times L_{\G}^2[0,T]\times L^2(\mu)$ to the BSDE 
\beqs
\tilde Y_t & = & \exp(\alpha B)+\int_t^T\tilde f(s, \tilde Y_s, \tilde Z_s, \tilde U_s)ds-\int_t^T \tilde Z_sdW_s-\int_t^T\int_E \tilde U_s(e)\mu(de,ds)\;,
\enqs
where the generator $\tilde f$ is defined by 
\beqs
\tilde f (t,y,z,u) & = & \inf_{\pi\in\Cc} \Big\{ \frac{\alpha^2}{2}|\pi\sigma_t|^2y- \alpha \pi\sigma_t (z+\vartheta_ty)+\int_E\Big[ e^{-\alpha\pi\beta_t(e)}( u(e)+y)-y\Big]\lambda_t(e) de \Big\} \;.
\enqs

We then notice that
 
\vspace{2mm}

\ni$\bullet$ $\tilde f $ satisfies (HUQ1) since it is an infimum of linear functions in the variable $z$,  
 
\vspace{2mm}

\ni$\bullet$  $\tilde f$ satisfies (HUQ2). Indeed, from the definition of $\tilde f $ we have
\beqs
\tilde f (t,y,z,u(.)-y)- \tilde f (t,y',z,u(.)-y') & \geq & 
\inf_{\pi\in\Cc}\Big\{(y-y')(\vartheta_t+\frac{\alpha}{2}\pi\sigma_t)\alpha\pi\sigma_t\Big\} 
-(y-y')\int_E\lambda_t(e)de\;,  
\enqs
for all $(t,z,u)\in[0,T]\times \R\times Bor(E,\R)$ and $y,y'\in\R$.
Since $\Cc$ is compact, we get from (HBI) the existence of a constant $C$ such that
\beqs
\tilde f (t,y,z,u-y)- \tilde f (t,y',z,u-y') & \geq & -C|y-y'|\;.
\enqs   
Inverting $y$ and $y'$ we get the result. 
 
\vspace{2mm}

\ni$\bullet$  $\tilde f $ satisfies (HUQ3). Indeed, since $0\in\Cc$, we get from (HBI) the existence of a constant $C$ such that
\beqs
\tilde f(t,y,z,u) & \leq & C\Big(|y|+\int_E|u(e)|\lambda_t(e)de\Big)\;, \quad (t,y,z,u)\in[0,T]\times\R\times\R\times Bor(E,\R)\;.
\enqs
We get from (HBI), there exists a positive constant $C$ s.t.
\beqs
\tilde f(t,y,z,u) & \geq & \inf_{\pi\in\Cc} \Big\{ \frac{\alpha^2}{2}|\pi\sigma_t|^2y - \alpha \pi\sigma_t ( z+\vartheta_t y)\Big\}\\
 & & +\inf_{\pi\in\Cc}\Big\{\int_Ee^{-\alpha \pi \beta_t(e)}(u(e)+y)\lambda_t(e)de\Big\}- C|y| \;.
\enqs
Then, from (HS1), (HS2) and the compactness of $\Cc$, we get
\beqs
\tilde f(t,y,z,u) & \geq & -C\Big(1+|y|+|z|+\int_E|u(e)|\lambda_t(e)de\Big)\;,\quad (t,y,z,u)\in[0,T]\times\R\times\R\times Bor(E,\R)\;.
\enqs
 
\vspace{2mm}

\ni$\bullet$  $\tilde f $ satisfies (HUQ4) since at time $t$  it is an integral of the variable $u$ w.r.t. $\lambda_t$, which vanishes on the interval $(\tau_n,\infty)$.   
 
\vspace{2mm}

Since $\tilde f $ satisfies (HUQ1), (HUQ2), (HUQ3) and (HUQ4), we get from Theorem \ref{TUQ} that $(\tilde Y^1,\tilde Z^1,\tilde U^1)=(\tilde Y^2,\tilde Z^2,\tilde U^2)$ in $\Sc^\infty_\G[0,T]\times L^2_\G[0,T]\times L^2(\mu)$. From the definition of $(\tilde Y^i,\tilde Z^i,\tilde U^i)$ for $i=1,2$, we get $(Y^1,Z^1,U^1)=(Y^2,Z^2,U^2)$ in $\Sc^\infty_\G[0,T]\times L^2_\G[0,T]\times L^2(\mu)$. 


\vspace{2mm}

\ni\textbf{Step 3.} We check that $M^{(\hat \pi)}$ is a BMO-martingale.  Since $\Cc$ is compact, (HS1) holds and $U$ is bounded as the jump part of the bounded process $Y$, it suffices to prove that $\int_0^.Z_sdW_s$ is a BMO-martingale. 

Let $M$ denote the upper bound of the uniformly bounded process $Y$.  
Applying It\^{o}'s formula to $(Y - M)^2$, we obtain for any stopping time $\tau \leq T$
\beqs
\E \Big[ \int_\tau^T |Z_s| ^2 ds \Big| \Gc_\tau \Big] & = & \E \big[ ( \xi - M )^2 \big| \Gc_\tau \big] - | Y_\tau - M |^2 \\
& & + 2 \E \Big[ \int_\tau^T (Y_s - M) f(s,  Z_s,U_s) ds \Big| \Gc_\tau \Big] \;.
\enqs
The definition of $f$ yields 
\beqs
- \vartheta_t Z_t - \frac{|\vartheta_t|^2}{2 \alpha} - \frac{1}{\alpha} \int_E \lambda_t(e) de & \leq & f(t,Z_t, U_t) \;,
\enqs
 for all $t \in [0, T]$. 
Therefore, since (HBI) and (HS4) hold, we get
\beqs
\E \Big[ \int_\tau^T |Z_s| ^2 ds \Big| \Gc_\tau \Big] & \leq & C\Big(1 +  \E \Big[ \int_\tau^T |Z_s + 1 | ds \Big| \Gc_\tau \Big]\Big)  \\
&\leq & C + \frac{1}{2}\E \Big[ \int_\tau^T |Z_s| ^2 ds \Big| \Gc_\tau \Big] \;.
\enqs
Hence, $\int_0^. Z_s dW_s$ is a BMO-martingale for $k= 0, \ldots, n$.

\vspace{2mm}

\ni\textbf{Step 4.} It remains to show that $R^{(\pi)}$ is a supermartingale for any $\pi\in\Ac$. Since $\pi\in\Ac$, the process $\mathfrak{E}(M^{(\pi)})$ is a positive local martingale, because it is the Doleans-Dade exponential of a local martingale whose the jumps are grower than $-1$. Hence, there exists a sequence of stopping times $(\delta_n)_{n\in\N}$ satisfying $\lim_{n\rightarrow \infty}\delta_n=T,~\P-a.s.$, such that $\mathfrak{E}(M^{(\pi)})_{.\wedge \delta_n}$ is a positive martingale for each $n\in \N$. The process $A^{(\pi)}$ is nondecreasing. Thus, $R^{(\pi)}_{t\wedge \delta_n}=R_0\mathfrak{E}(M^{(\pi)})_{t\wedge \delta_n}\exp(A^{(\pi)}_{t\wedge\delta_n})$ is a supermartingale, i.e. for $s\leq t$
\beqs
\E \big[R^{(\pi)}_{t\wedge \delta_n} \big|\Gc_s \big] \leq R^{(\pi)}_{s\wedge \delta_n}\;.
\enqs
For any set $A\in\Gc_s$, we have 
\beq \label{3inequation R}
\E \big[R^{(\pi)}_{t\wedge \delta_n}\mathds{1}_A \big] \leq \E \big[R^{(\pi)}_{s\wedge \delta_n}\mathds{1}_A \big]\;.
\enq
On the other hand, since
\beqs
R^{(\pi)}_t = -\exp\big(-\alpha (X^{x, \pi}_t - Y_t)\big)\;,
\enqs
we use both the uniform integrability of $(\exp(-\alpha X^{x,\pi}_\delta))$ where $\delta$ runs over the set of all stopping times and the boundedness of $Y$ to obtain the uniform integrability of 
\beqs
\{R^{(\pi)}_{\tau}~:~\tau \text{ stopping time valued in } [0,T]\}.
\enqs
 Hence, the passage to the limit as $n$ goes to $\infty$ in \reff{3inequation R} is justified and it implies
\beqs
\E\big[R^{(\pi)}_{t}\mathds{1}_A \big] \leq \E\big[R^{(\pi)}_{s}\mathds{1}_A \big] \;.
\enqs
We obtain the supermartingale property of $R^{(\pi)}$.\\
To complete the proof, we show that the strategy $\hat{\pi}$ defined by \reff{pi optimal} is optimal. We first notice that from Lemma \ref{lem ttes adm} we have $\hat \pi\in\Ac$. By definition of $\hat{\pi}$, we have $A^{(\hat{\pi})} = 0$ and hence, $R^{(\hat{\pi})}_t = R_0\mathfrak{E}(M^{(\hat{\pi})})_t$. 
Since $\Cc$ is compact,  (HS1) holds and $U$ is bounded as jump part of the bounded process $Y$, there exists a constant $\delta>0$ s.t. 
\beqs
\Delta M^{(\hat \pi)}_t & = & M^{(\hat \pi)}_t-M^{(\hat \pi)}_{t^-}\geq -1+\delta\;. 
\enqs 
Applying  Kazamaki criterion to the BMO martingale $M^{(\hat \pi)}$ (see \cite{kaz79}) we obtain that $\mathfrak{E}(M^{(\hat \pi)})$  
 is a true martingale. As a result, we get
\beqs
\sup_{\pi\in \Ac} \E \big(R^{(\pi)}_T\big)  =  R_0 = V(x)\;.
\enqs
Using that $(Y, Z, U)$ is the unique solution of the BSDE  \reff{3equation exponentielle}, we obtain the expression (\ref{3fonction valeur}) for the value function.
\ep

\begin{Remark}
{ \rm Concerning the existence and uniqueness of a solution to BSDE \reff{3equation exponentielle},  we notice that the compactness assumption on $\Cc$ is only need for the uniqueness. Indeed, in the case where $\Cc$ is only a closed set, the generator of the BSDE still satisfies a quadratic growth condition which allows to apply Kobylanski existence result.  However, for the uniqueness of the solution to BSDE \reff{3equation exponentielle}, we need $\Cc$ to be compact to get  Lipschitz continuous decomposed generators w.r.t. $y$.  We notice that the existence result for a similar BSDE in the case of Poisson jumps is proved by Morlais in \cite{mor09} and \cite{mor10} without any compactness assumption on $\Cc$. }
\end{Remark}
\appendix

\vspace{2cm}

\ni\textbf{{\huge Appendix}}
\section{Proof of Lemma \ref{3lemme decomposition} (ii)}
\setcounter{equation}{0} \setcounter{Assumption}{0} 
\setcounter{Theorem}{0} \setcounter{Proposition}{0}
\setcounter{Corollary}{0} \setcounter{Lemma}{0}
\setcounter{Definition}{0} \setcounter{Remark}{0}
We prove the decomposition for the progressively measurable processes $X$ of the form
\beqs
X_t & = & J_t+\int_0^tU_s(e)\mu(de,ds)\;,\quad t\geq 0\;,
\enqs
where $J$ is $\Pc(\G)$-measurable and $U$ is $\Pc(\G)\otimes\Bc(E)$-measurable. To prove the decomposition \reff{3decomposition progressive}, it sufficies to prove it for the process $J$ and the process $V$ defined by
\beqs
V_t & = & \int_0^tU_s(e)\mu(de,ds)\;,\quad t\geq 0\;.
\enqs
$\bullet$ Decomposition of the process $J$.

\ni Since $J$ is $\Pc(\G)$-measurable, we can write
\beqs
J_t & = & J^0_t\mathds{1}_{t\leq \tau_1}+\sum_{k=1}^{n}J^k_t( \tau_{(k)},\zeta_{(k)})\mathds{1}_{\tau_k<t\leq \tau_{k+1}}\;,
\enqs
for all $t \geq 0$, where $J^0$ is $\mathcal{P}(\F)$-measurable, and $J^k$ is $\mathcal{P}(\F) \otimes\Bc(\Delta_{k})\otimes\Bc(E^k)$-measurable, for $k=1,\ldots,n$.  This leads to the following decomposition of $J$:
\beqs
J_t & = & J^0_t\mathds{1}_{t\leq \tau_1}+\sum_{k=1}^{n}\bar J^k_t( \tau_{(k)},\zeta_{(k)})\mathds{1}_{\tau_k\leq t< \tau_{k+1}}\;,
\enqs
where
\beqs
\bar J^k_t( \theta_{(k)},e_{(k)}) & = & J^k_t( \theta_{(k)},e_{(k)}) + \big( J^{k-1}_t( \theta_{(k-1)},e_{(k-1)})- J^k_t( \theta_{(k)},e_{(k)}) \big)\mathds{1}_{t=\theta_k}\;,
\enqs
for $k=1,\ldots,n$ and $(\theta_{(k)},e_{(k)})\in\Delta_k\times E^k$. Since $J^k$ is $\Pc(\F)\otimes\Bc(\Delta_k)\otimes\Bc(E^k)$-measurable for all $k=0,\ldots,n$, we get that $(\bar J^k_t)_{t\in[0,s]}$ is $\Fc_s\otimes \Bc([0,s])\otimes\Bc(\Delta_k)\otimes\Bc(E^k)$-measurable for all $s\geq 0$.

\vspace{2mm}

\ni$\bullet$ Decomposition of the process $V$.

\ni Since $U$ is $\Pc(\G)\otimes\Bc(E)$-measurable, we can write
\beqs
U_t(.) & = & U^0_t(.)\mathds{1}_{t\leq \tau_1}+\sum_{k=1}^{n}U^k_t( \tau_{(k)},\zeta_{(k)},.)\mathds{1}_{\tau_k<t\leq \tau_{k+1}}\;,
\enqs
for all $t \geq 0$, where $U^0$ is $\mathcal{P}(\F)\otimes\Bc(E)$-measurable, and $U^k$ is $\mathcal{P}(\F) \otimes\Bc(\Delta_{k})\otimes\Bc(E^k)\otimes\Bc(E)$-measurable, for $k=1,\ldots,n$.  This leads to the following decomposition of $V$:
\beqs
V_t & = & \sum_{k=1}^{n}U^{k-1}_{\tau_{k}}(\tau_{(k-1)},\zeta_{(k)})\mathds{1}_{\tau_{k}\leq t}\\
 & = & \sum_{k=1}^{n}\Big(\sum_{j=1}^kU^{j-1}_{\tau_{j}}(\tau_{(j-1)},\zeta_{(j)})\mathds{1}_{\tau_{j}\leq t}\Big)\mathds{1}_{\tau_{k}\leq t<\tau_{k+1}}\\
  & = &\sum_{k=1}^{n}V^k_{t}(\tau_{(k)},\zeta_{(k)})\mathds{1}_{\tau_{k}\leq t<\tau_{k+1}} \;,
\enqs
where $V^k$ is defined by $V^0=0$ and
\beqs
V^k_{t}(\theta_{(k)},e_{(k)})  & = & \sum_{j=1}^kU^{j-1}_{\theta_{j}}(\theta_{(j-1)},e_{(j)})\mathds{1}_{\theta_{j}\leq t}\;,\quad t\geq 0\;,~(\theta_{(k)},e_{(k)})\in\Delta_k\times E^k\;,
\enqs
for $k=1,\ldots,n$.
We now check that for all $s\geq 0$, $(V^k_t(.))_{t\in[0,s]}$ is $\Fc_s\otimes\Bc([0,s])\otimes\Bc(\Delta_k)\otimes\Bc(E^k)$-measurable. Since $U^j$ is $\Pc(\F)\otimes \Bc(\Delta_j)\otimes\Bc(E^j)$-measurable, we get that $(U^j_t(.))_{t\in[0,s]}$ is $\Fc_s\otimes\Bc([0,s])\otimes\Bc(\Delta_j)\otimes\Bc(E^j)$-measurable. Therefore $(t,\theta_{(j)},e_{(j)})\in[0,s]\times\Delta_j\times E^j\mapsto U^{j-1}_{\theta_j}(\theta_{(j-1)},e_{(j)})\mathds{1}_{\theta_{j}\leq t}$ is $\Fc_s\otimes\Bc([0,s])\otimes\Bc(\Delta_j)\otimes\Bc(E^j)$ for $j =0,\ldots,n$. From the definition of $V^k$ we get that $(V^k_t(.))_{t\in[0,s]}$ is $\Fc_s\otimes\Bc([0,s])\otimes\Bc(\Delta_k)\otimes\Bc(E^k)$-measurable.
\ep

\section{Proof of Proposition \ref{prop intensite}} 
 \setcounter{equation}{0} \setcounter{Assumption}{0} 
\setcounter{Theorem}{0} \setcounter{Proposition}{0}
\setcounter{Corollary}{0} \setcounter{Lemma}{0}
\setcounter{Definition}{0} \setcounter{Remark}{0}

We first give a lemma which is a generalization of a proposition  in \cite{nekjeajia}. Throughout the sequel, we denote
\beqs
\Ec_t^{\F,i,k}\big(G\big)(\theta_{(i-1)},e_{(i-1)}) &= & \int_{\Delta_{k-i+1}\times E^{k-i+1}}\mathds{1}_{\theta_i>t} \E \big[G ( \theta_{(k)}, e_{(k)}) \big| \Fc_t \big]d \theta_{i}\ldots d \theta_{k}  d e_{i}\ldots d e_{k} 
\;,
\enqs
for any $\Fc_\infty\otimes\Bc(\Delta_k)\otimes\Bc(E^k)$-measurable function $G$ and any integers $i$ and $k$ such that $1\leq i\leq k\leq n$.

\begin{Lemma}\label{corollaire}
Fix $t,s\in\R_+$ such that $t\leq s$. Let $X$ be a positive $\Fc_s \otimes \Bc ( \Delta_n) \otimes \Bc (E^n)$-measurable function on $\Omega \times \Delta_n \times E^n$, then
\beqs
\E \big[X(\tau_{(n)}, \zeta_{(n)}) \big| \Gc_t \big] & = & \sum_{i = 0}^{n} \mathds{1}_{\tau_{i} \leq t < \tau_{i+1}}\frac{\Ec_t^{\F,i+1,n}\big(X\gamma_s\big)(\tau_{(i)},\zeta_{(i)})}{\Ec_t^{\F,i+1,n}\big(\gamma_t\big)(\tau_{(i)},\zeta_{(i)})}\;.\\
\enqs

\end{Lemma}

\noindent\textbf{Proof.}
Let $H$ be a positive and $\Gc_t$-measurable test random variable, which can be written
\beqs
H & = & \sum_{i=0}^nH^i(\tau_{(i)}, \zeta_{(i)})\mathds{1}_{ \tau_i \leq t < \tau_{i+1} } \;,
\enqs
where $H^i$ is $\Fc_t \otimes \Bc(\Delta_i) \otimes \Bc(E^i)$-measurable for $i=0,\ldots,n$. Using the joint density $\gamma_t(\theta, e)$ of $(\tau, \zeta)$, we have on the one hand
\beqs
 \E [ \mathds{1}_{ \tau_i \leq t < \tau_{i+1}} H X(\tau_{(n)}, \zeta_{(n)})]  
&  = &  \E \Big[ \int_{(0,t]^i\cap\Delta_i\times E^i}  \hspace{-5mm}d\theta_{(i)} de_{(i)} H^i_t(\theta_{(i)}, e_{(i)}) \Ec_t^{\F,i+1,n}\big(X\gamma_s\big)(\tau_{(i)},\zeta_{(i)})\Big]\;.
\enqs
 On the other hand, we have 
\beqs
&  & \E \Big[ \mathds{1}_{ \tau_i \leq t < \tau_{i+1}} H \frac{\Ec_t^{\F,i+1,n}\big(X\gamma_s\big)(\tau_{(i)},\zeta_{(i)})}{\Ec_t^{\F,i+1,n}\big(\gamma_t\big)(\tau_{(i)},\zeta_{(i)})} \Big]\\
& = & \E \Big[ \mathds{1}_{ \tau_i \leq t < \tau_{i+1}} H^i(\tau_{(i)},\zeta_{(i)}) \frac{\Ec_t^{\F,i+1,n}\big(X\gamma_s\big)(\tau_{(i)},\zeta_{(i)})}{\Ec_t^{\F,i+1,n}\big(\gamma_t\big)(\tau_{(i)},\zeta_{(i)})} \Big]\\
& =  & \E \Big[ \int_{(0,t]^i\cap\Delta_i\times E^i}  \hspace{-5mm}d\theta_{(i)} de_{(i)} H^i_t(\theta_{(i)}, e_{(i)}) \frac{\Ec_t^{\F,i+1,n}\big(X\gamma_s\big)(\theta_{(i)}, e_{(i)})}{\Ec_t^{\F,i+1,n}\big(\gamma_t\big)(\theta_{(i)}, e_{(i)})} 
\Ec_t^{\F,i+1,n}\big(\gamma_t\big)(\theta_{(i)}, e_{(i)})\Big]\\
&= &  \E [ \mathds{1}_{ \tau_i \leq t < \tau_{i+1}} H X(\tau_{(n)}, \zeta_{(n)})]  \;.
\enqs
\ep\\

We now prove Proposition \ref{prop intensite}. To this end, we prove that for any nonnegative $\Pc(\G)\otimes\Bc(E)$-measurable process $U$, any $T>0$   and any $t\in[0,T]$, we have
\beq\label{eg a montrer}
\E \Big[ \int_t^T \int_E U_s(e)\mu(de, ds) \Big| \Gc_t \Big] & = & \E \Big[ \int_t^T \int_E U_s(e)\lambda_s (e) de ds \Big| \Gc_t \Big] \;,
\enq
where $\lambda$ is defined by \reff{decomp lambda}.
 
We first study the left hand side of \reff{eg a montrer}. From Lemma \ref{3lemme decomposition} and Remark \ref{remparam}, we can write 
\beqs
U_t(e) & = & \sum_{k=0}^n\mathds{1}_{\tau_k< t \leq \tau_{k+1}}U^k_t(\tau_{(k)},\zeta_{(k)},e)\;,\quad (t,e)\in [0,T]\times E\;, 
\enqs
where $U^k$ is a $\Pc(\G)\otimes\Bc(\Delta_k)\otimes\Bc(E^{k+1})$-measurable process for $k=0,\ldots,n$. Moreover, since $U$ is nonnegative, we can assume that $U^k$, $k=0,\ldots,n$, are nonnegative.
Then, from Lemma \ref{corollaire}, we have:
\beqs
 & & \E \Big[ \int_t^T \int_E U_s(e)\mu(de, ds) \Big| \Gc_t \Big] 
~=~
\sum_{k = 1}^n \E \big[ \mathds{1}_{t < \tau_k \leq T}U^{k-1}_{\tau_k}(\tau_{(k-1)},\zeta_{(k)}) \big| \Gc_t \big]  \\
 & = & \sum_{k = 1}^n  \sum_{i = 0}^n \mathds{1}_{\tau_i \leq t < \tau_{i+1}}  \frac{ \Ec^{\F,i+1, n}_t \Big(\mathds{1}_{t < \theta_k \leq T}U^{k-1}_{\theta_k}(\theta_{(k-1)},e_{(k)}) \gamma_T( \theta, e)\Big)(\tau_{(i)},e_{(i)})}{\Ec^{\F,i+1,n}_t \big(\gamma_t\big)(\tau_{(i)},e_{(i)})} \\ 
 & = & \sum_{k, i = 0\atop  i\leq k}^{n-1} 
\mathds{1}_{\tau_i \leq t < \tau_{i+1}}  \frac{\Ec^{\F,i+1,n}_t \Big(\mathds{1}_{t < \theta_{k+1} \leq T}U^{k}_{\theta_{k+1}}(\theta_{(k)},e_{(k+1)}) \gamma_{\theta_{k+1}}( \theta, e)\Big)(\tau_{(i)},e_{(i)})}{\Ec^{\F,i+1,n}_t \big(\gamma_t\big)(\tau_{(i)},e_{(i)})}   \\
&=  & \sum_{k, i = 0\atop  i\leq k}^{n-1}  
\mathds{1}_{\tau_i \leq t < \tau_{i+1}} \frac{
\Ec^{\F,i+1,k+1}\Big(
\mathds{1}_{t < \theta_{k+1}\leq T}U^{k}_{\theta_{k+1}}(\theta_{(k)},e_{(k+1)}) \gamma^{k+1}_{\theta_{k+1}}( \theta_{(k+1)}, e_{(k+1)})  \Big)(\tau_{(i)},e_{(i)})
 }{
\Ec^{\F,i+1,n}\big(\gamma_t\big)(\tau_{(i)},e_{(i)}) } \;.
\enqs

\ni We now study the right hand side of \reff{eg a montrer}:
\beqs
&  &\E \Big[ \int_t^T \int_E U_s(e)\lambda_s (e) de ds \Big| \Gc_t \Big]  =   \sum_{k=0}^{n-1} \E \Big[ \int_t^T\hspace{-3mm} \int_E\hspace{-1mm} \mathds{1}_{\tau_{k} < s \leq \tau_{k+1}} U^k_s( \tau_{(k)}, \zeta_{(k)})\lambda^{k+1}_s ( e, \tau_{(k)}, \zeta_{(k)}) de ds \Big| \Gc_t \Big]\\
  & = & \sum_{k= 0}^{n-1}  \sum_{i = 0}^n \mathds{1}_{\tau_i \leq t < \tau_{i+1}} \hspace{-1mm}\frac{ \Ec_t^{\F,i+1,n}\Big(\hspace{-1mm} \int_t^T\hspace{-1mm} \int_E\hspace{-0.5mm}\mathds{1}_{ \theta_{k} < s \leq \theta_{k+1}} U^k_s(\theta_{(k)} , e_{(k)})\lambda^{k+1}_s ( e', \theta_{(k)} , e_{(k)}) \gamma_s( \theta, e)de' ds \Big) (\tau_{(i)},\zeta_{(i)})}{\Ec^{\F,i+1,n}_t \big(\gamma_t\big)(\tau_{(i)},e_{(i)})} \\
 & = & \sum_{k,i= 0 \atop i\leq k}^{n-1} 
 \mathds{1}_{\tau_i \leq t < \tau_{i+1}} \hspace{-1mm}\frac{ \Ec_t^{\F,i+1,k}\Big( \int_t^T\hspace{-1mm} \int_E  \mathds{1}_{ \theta_{k} < s } U^k_s(\theta_{(k)} , e_{(k)})\lambda^{k+1}_s ( e, \theta_{(k)} , e_{(k)}) \gamma_s^k( \theta_{(k)}, e_{(k)})  de' ds\Big)(\tau_{(i)},\zeta_{(i)})}{\Ec^{\F,i+1,n}_t \big(\gamma_t\big)(\tau_{(i)},e_{(i)})} \\
 & =  & \sum_{k,i= 0 \atop i\leq k}^{n-1} 
 \mathds{1}_{\tau_i \leq t < \tau_{i+1}} \frac{ \Ec_t^{\F,i+1,k}\Big( \int_t^T \hspace{-1mm}\int_E  \mathds{1}_{ \theta_{k} < s } U^k_s(\theta_{(k)} , e_{(k)})\gamma^{k+1}_s (\theta_{(k)},s , e_{(k)},e') de' ds\Big) }{\Ec^{\F,i+1,n}_t \big(\gamma_t\big)(\tau_{(i)},e_{(i)})}\;, 
\enqs
where the last equality comes from the definition of $\lambda^k$. Hence, we get \reff{eg a montrer}.

\section{Measurability of solutions to BSDEs depending on a parameter}

\subsection{Representation for Brownian martingale depending on a parameter}
\setcounter{equation}{0} \setcounter{Assumption}{0} 
\setcounter{Theorem}{0} \setcounter{Proposition}{0}
\setcounter{Corollary}{0} \setcounter{Lemma}{0}
\setcounter{Definition}{0} \setcounter{Remark}{0}

We consider $\Xc$ a Borelian subset of $\R^p$ and $\rho$ a finite measure on $\Bc(\Xc)$. Let  $\{\xi(x)~:~x\in \Xc\}$ be a family of random variables such that the map $\xi:~\Omega\times\Xc\rightarrow \R$ is   $\Fc_T\otimes \Bc(\Xc)-$measurable and satisfies $\int_\Xc\E|\xi(x)|^2\rho(dx)<\infty$. In the following result, we generalize the representation property as a stochastic integral w.r.t.  $W$ of square-integrable random variables to the family $\{\xi(x)~:~x\in \Xc\}$. The proof follows the same lines as for the classical It\^o representation Theorem which can be found e.g. in \cite{Ok07}. For the sake of completeness we sketch the proof. 
\begin{Theorem}\label{repxi}
There exists a $\Pc(\F)\otimes \Bc(\Xc)$-measurable map $Z$ such that $\int_\Xc\int_0^T\E|Z_s(x)|^2ds\rho(dx)$ $<$ $\infty$ and 
\beq\label{repxiZ}
\xi(x) & = & \E[\xi(x)] + \int_0^T Z_s(x)dW_s\;, \quad \P\otimes\rho-a.e. 
\enq
\end{Theorem}
As for the standard representation theorem, we first need a lemma which provides a dense subset of $L^2(\Fc_T\otimes\Bc(\Xc), \P\otimes\rho)$ generated by easy functions.
\begin{Lemma}
Random variables of the form
\beq\label{ssensdense}
\exp\Big(\int_0^Th_t(x)dW_t-{1\over2} \int_0^T|h_t(x)|^2dt\Big)\;,
\enq
where $h$ is a bounded  $\Bc([0,T])\otimes\Bc(\Xc)-$measurable map span a dense subset of $L^2(\Fc_T\otimes\Bc(\Xc), \P\otimes\rho)$.
\end{Lemma}
\noindent\textbf{Sketch of the proof.} Let $\Lambda\in L^2(\Fc_T\otimes\Bc(\Xc), \P\otimes\rho)$ orthogonal to all functions of the form \reff{ssensdense}. Then, in particular, we have
\beqs
G(\alpha_1,\ldots,\alpha_n)~=~\int_\Xc\E\big[\Lambda\exp(\alpha_1W_{t_1}+\cdots+\alpha_nW_{t_n})\big]d\rho & = & 0 \;,
\enqs
for all $\alpha_1,\ldots,\alpha_n\in\R$ and all $t_1,\ldots,t_n\in[0,T]$.  Since $G$ is identically equal to zero on $\R^n$ and is analytical it is also identically equal to $0$ on $\mathbb{C}^n$. We then have for any $\Bc(\Xc)\otimes\Bc(\R^p)-$ measurable function $\phi$ such that $\phi(x,.)\in C^\infty(\R^n)$ with compact support for all $x\in\Xc$ 
\beqs
\int_\Xc\E[Y\phi(x,W_{t_1},\ldots,W_{t_n})]d\rho (x)& = & \\
\int_{\R^n\times\Xc}\hat \phi(x,\alpha_1,\ldots,\alpha_n)\E\big[\Lambda\exp(\alpha_1W_{t_1}+\cdots+\alpha_nW_{t_n})\big]d\rho(x)d\alpha_1\ldots d\alpha_n & =& 0 \;,
\enqs
where $\hat\phi(x,.)$ is the Fourier transform of $\phi(x,.)$.  Hence, $\Lambda$ is equal to zero since it is orthogonal to a dense subset of $L^2(\Fc_T\otimes \Bc(\Xc))$. 
\ep

\vspace{2mm}

\noindent\textbf{Sketch of the proof of Theorem \ref{repxi}.} First suppose that $\xi$ has the following form:
\beqs
\xi(x) & = & \exp\Big(\int_0^Th_t(x)dW_t-{1\over 2}\int_0^T|h_t(x)|^2dt\Big)\;,
\enqs 
with $h$ a bounded $\Bc([0,T])\otimes\Bc(\Xc)-$measurable map. Then, applying It\^o's formula to the process $\exp\Big(\int_0^.h_t(x)dW_t-{1\over 2}\int_0^.|h_t(x)|^2dt\Big)$, we get that $\xi$ satisfies \reff{repxiZ}
where the process $Z$ is given by
\beqs
Z_t(x) & = & h_t(x)\exp\Big(\int_0^th_s(x)dW_s-{1\over 2}\int_0^t|h_s(x)|^2ds\Big)\;, \quad (t,x)\in[0,T]\times\Xc\;.
\enqs
Now for any $\xi\in L^2(\Fc_T\otimes \Bc(\Xc),\P\otimes\rho)$, there exists a sequence $(\xi^n)_{n\in \N}$ such that 
each $\xi^n$ satisfies 
\beqs
\xi^n(x) & = & \E[\xi^n(x)] + \int_0^T Z^n_s(x)dW_s\;, \quad \P\otimes\rho-a.e. 
\enqs
and $(\xi^n)_{n\in \N}$ converges to $\xi$ in $L^2(\Fc_T\otimes \Bc(\Xc),\P\otimes dt \otimes\rho)$ . Then, using It\^o's Isometry, we get that the sequence $(Z^n)_{n\in \N}$ is Cauchy and hence converges in $L^2(\Pc(\F)\otimes \Bc(\Xc),\P\otimes dt \otimes\rho)$ to some $Z$. Using again the It\^o Isometry, we get that 
 $(\xi^n)_{n \in \N}$ converges to $\E[\xi(x)]+\int_0^TZ_s(x)dW_s$ in $L^2(\Fc_T\otimes \Bc(\Xc),\P\otimes\rho)$. Identifying the limits, we get the result. 
\ep
\begin{Corollary}\label{corMG}
Let $M$ be a $\Pc(\F)\otimes \Bc(\Xc)-$measurable map such that $(M_t(x))_{0\leq t\leq T}$ is a martingale for all $x\in\Xc$ and $\int_\Xc\E |M_T(x)|^2\rho(dx)<\infty$. Then, there exists a $\Pc(\F)\otimes\Bc(\Xc)-$measurable map $Z$ such that $\int_0^T \int_\Xc \E|Z_s(x)|^2 \rho(dx) ds$ $<$ $\infty$ and  
\beqs
M_t(x) & = & M_0(x) +\int_0^tZ_s(x)dW_s\;.
\enqs
\end{Corollary}
\noindent The proof is a direct consequence of Theorem \ref{repxi} as in \cite{Ok07} so we omit it. 
\subsection{BSDEs depending on a parameter}
We now study the measurability of solutions to Brownian BSDEs whose data depend on the  parameter $x\in\Xc$. We consider 
\begin{itemize}
\item a family  $\{\xi(x)~:~x\in \Xc\}$ of random variables such that the map $\xi:~\Omega\times\Xc\rightarrow \R$ is   $\Fc_T\otimes \Bc(\Xc)-$measurable and satisfies $\int_\Xc\E|\xi(x)|^2\rho(dx)<\infty$,

\item  a family  $\{f(., x)~:~x\in \Xc\}$ of random maps such that the map $f:~\Omega\times[0,T]\times\R\times\R^d\times\Xc\rightarrow \R$ is $\Pc(\F)\otimes\Bc(\R)\otimes \Bc(\R^d)\otimes \Bc(\Xc)-$measurable and satisfies $\int_0^T \int_\Xc \E|f(s,0,0,x)|^2 \rho(dx) ds<\infty$.
\end{itemize}

We then consider the BSDEs depending on the parameter $x\in\Xc$:
\beq\label{EDSR param}
Y_t(x) & = & \xi(x)+\int_t^Tf(s,Y_s(x),Z_s(x),x)ds-\int_t^TZ_s(x)dW_s \;, \quad  (t,x)\in [0,T]\times\Xc\;.\qquad
\enq
\begin{Lemma}\label{Lem Pic}
Assume that the generator $f$ does not depend on $(y,z)$ i.e. $f(t,y,z,x)=f(t,x)$. Then, BSDE \reff{EDSR param} admits a solution $(Y,Z)$ such that $Y$ and $Z$ are $\Pc(\F)\otimes\Bc(\Xc)-$measurable.   
\end{Lemma}
\noindent\textbf{Proof.}
Consider the family of martingales $\{M(x)~:~x\in\Xc\}$, where $M$ is defined by
\beqs
M_t(x) & = & \E\Big[\xi(x)+\int_0^Tf(s,x)ds\Big|\Fc_t\Big]\;, \quad  (t,x)\in [0,T]\times\Xc\;.
\enqs
Then, from Corollary \ref{corMG}, there exists a $\Pc(\F)\otimes \Bc(\R^d)-$measurable map $Z$ such that 
$\int_0^T\int_\Xc \E|Z_s(x)|^2 \rho(dx) ds$ $<$ $\infty$ and  
\beqs
M_t(x) & = & M_0(x) +\int_0^tZ_s(x)dW_s\;, \quad (t,x)\in [0,T]\times\Xc\;.
\enqs
We then easily check that the process $Y$ defined by 
\beqs
Y_t(x) & = & M_t(x)-\int_0^tf(s,x)ds\;, \quad (t,x) \in [0,T]\times\Xc\;,
\enqs
is $\Pc(\F)\otimes\Bc(\Xc)-$measurable and that $(Y,Z)$ satisfies \reff{EDSR param}. 
\ep

\vspace{2mm}

We now consider the case where the generator $f$ is Lipschitz continuous: there exists a constant $L$ such that 
\beq\label{condLip}
|f(t,y,z,x)-f(t,y',z',x)| & \leq & L(|y-y'|+|z-z'|)\;,
\enq
for all  $(t,y,y',z,z')\in[0,T]\times[\R]^2\times[\R^d]^2$.

\begin{Proposition}\label{mesLipsol}
Suppose that $f$ satisfies \reff{condLip}. Then, BSDE \reff{EDSR param} admits a $\Pc(\F)\otimes\Bc(\Xc)-$measurable solution $(Y,Z)$ such that $\E\int_0^T\int_\Xc(|Y_s(x)|^2+|Z_s(x)|^2) \rho(dx) ds<\infty$. 
\end{Proposition}
\ni \textbf{Proof.} 
 Consider the sequence $(Y^n,Z^n)_{n\in \N}$ defined by $(Y^0,Z^0)= (0,0)$ and for $n\geq 1$
 \beqs
 Y^{n+1}_t(x) & =& \xi(x)+\int_t^Tf(s,Y^n_s(x),Z_s^n(x))ds-\int_t^TZ^{n+1}_s(x)dW_s,\; (t,x)\in[0,T]\times\Xc\;.
 \enqs
 From Lemma \ref{Lem Pic}, we get that $(Y^n,Z^n)$ is $\Pc(\F)\otimes\Bc(\Xc)-$measurable for all $n\in\N$.  Moreover, since  $f$ satisfies \reff{condLip}, the sequence $(Y^n,Z^n)_{n\in \N}$ converges (up to a subsequence) a.e. to $(Y,Z)$ solution to \reff{EDSR param} (see \cite{pp90}). Hence, the solution $(Y,Z)$ is also $\Pc(\F)\otimes\Bc(\Xc)-$measurable.\ep
 
 \section{A regularity result for the decomposition}
 
 \begin{Proposition}\label{reg-decomp}
 Let $p\geq1$ and $(f_t(x))_{(t,x)\in[0,T]\times \R^p}$ be a $\Pc(\G)\otimes\Bc(\R^p)$-measurable map. Suppose that $f_t(.)$ is locally uniformly continuous (uniformly in $\omega\in \Omega$). Then $f^k_t(.,\theta_{(k)},e_{(k)})$ is locally uniformly continuous (uniformly in $\omega\in \Omega$) for $\theta_k\leq t$ and $k=0,\ldots,n$. 
 \end{Proposition}
 \ni \textbf{Proof.} For sake of clarity, we prove the result without marks, but the argument easily extends to the case with marks. Fix $k\in\{0,\ldots,n\}$ and  for $R>0$,  denote by $mc_R^f$ the modulus of continuity of $f$ on $B_{\R^p}(0,R)$. 
 Then for any $\tilde \theta_k>\cdots>\tilde \theta_1>0$ and $h_1,\ldots, h_n>0$  we have 
 from the definition of $mc_R^f$ and (HD)
 \beqs
 {1\over h_1\cdots h_k}\E\Big[|f_t(x)-f_t(x')|\mathds{1}_{\cap_{\ell\leq k}\{\tilde \theta_\ell-h_\ell\leq \tau_\ell\leq \tilde \theta_\ell<t\leq \tau_{\ell+1}\}}\Big|\Fc_t\Big] & \leq &\\
  mc^f_R(\eps) {{1\over h_1\cdots h_k}}\int_{\tilde \theta _1-h_1}^{\tilde \theta _1}d\theta_1\ldots\int_{\tilde \theta _k-h_k}^{\tilde \theta _k}d\theta_k\Big(\int\gamma_t(\theta)d\theta_{k+1}\ldots d\theta_n\Big)	\;,
 \enqs
for $x,x'\in B_{\R^p}(0,R)$ s.t. $|x-x'|\leq \eps$. Using the decomposition of $f$ we have
 \beqs
{1\over h_1\cdots h_k}\E\Big[|f_t(x)-f_t(x')|\mathds{1}_{\cap_{\ell\leq k}\{\tilde \theta_\ell-h_\ell \leq \tau_\ell\leq \tilde \theta_\ell<t\leq \tau_{\ell+1}\}}\Big|\Fc_t\Big] & = &\\
{1\over h_1\cdots h_k}\int_{\tilde \theta _1-h}^{\tilde \theta _1}d\theta_1\ldots\int_{\tilde \theta_k-h}^{\tilde\theta_k}
|f_t^k(x,\theta_{(k)})-f_t^k(x',\theta_{(k)})|\Big(\int\gamma_t(\theta,e)d\theta_{k+1}\ldots d\theta_n \Big)	d\theta_k\;.
  \enqs
  Sending each $h_\ell$ to zero we get 
  \beqs
 |f_t^k(x,\tilde \theta_{(k)})-f_t^k(x',\tilde\theta_{(k)})| & \leq & mc^f_R(\eps)\;.
  \enqs
 \ep


\begin{thebibliography}{}
\bibitem{ankblaeyr09} Ankirchner S., Blanchet-Scalliet C. and A. Eyraud-Loisel (2009): ``Credit risk premia and quadratic BSDEs with a single jump", forthcoming in {\it International Journal of Theoretical and Applied Finance }. 


\bibitem{barlow} Barlow M. T. (1978): ``Study of a Filtration Expanded to Include an Honest Time'', {\it Z. Wahrscheinlichkeitstheorie verw. Gebiete}, {\bf 44}, 307-323. 


\bibitem{bierut04} Bielecki T. and M. Rutkowski (2004): ``Credit risk: modelling, valuation and hedging'', Springer Finance.

\bibitem{biejearut04a} Bielecki T., Jeanblanc M. and M. Rutkowski (2004): ``Stochastic Methods in Credit Risk Modelling'', Lectures notes in Mathematics, Springer, {\bf 1856}, 27-128.


\bibitem{biejea08} Bielecki T. and M. Jeanblanc (2008): ``Indifference prices in Indifference Pricing'', {\it Theory and Applications, Financial Engineering}, Princeton University Press. R. Carmona editor  volume.

\bibitem{bh06} Briand P. and Y. Hu (2006): ``BSDE with quadratic growth and unbounded terminal
value", {\it Probability Theory Related Fields}, {\bf 136 (4)}, 604-618.

\bibitem{bh08} Briand P. and Y. Hu (2008): ``Quadratic BSDEs with convex generators and unbounded terminal conditions'', {\it Probability Theory and Related Fields}, {\bf 141}, 543-567. 

\bibitem{gio10} Callegaro G. (2010): ``Credit risk models under partial information'', {\it PhD Thesis}  Scuola Normal Superior di Pisa and Université d'\'Evry Val d'Essonne.

\bibitem{delhuri11} Delbaen F., Hu Y. and A. Richou (2011): ``On the uniqueness of solutions to quadratic
BSDEs with convex generators and unbounded terminal conditions", {\it Ann.
Inst. Henri Poincar\'e Probab. Stat.}, {\bf 47 (2)}, 559-574.

\bibitem{demey75} Dellacherie C. and P.A. Meyer (1975): ``ProbabilitŽ\'{e}s et Potentiel - Chapitres I ˆ ˆ\`a IV'', Hermann, Paris. 

\bibitem{demey80} Dellacherie C. and P.A. Meyer (1980): ``ProbabilitŽ\'{e}s et Potentiel - Chapitres V ˆ\`a VIII'', Hermann, Paris. 

\bibitem{dufsin03} Duffie D. and K. Singleton (2003): ``Credit risk: pricing, measurement and management'', Princeton University Press. 


\bibitem{nekjeajia} El Karoui N., Jeanblanc M. and Y. Jiao (2010): ``Modelling Succesive Default Events'', {\it Preprint}.

\bibitem{HeW} He S., Wang J.  and  Yan J. (1992):  ``Semimartingale theory and stochastic calculus'', Science Press, CRC Press, New-York. 

\bibitem{huimkmul05} Hu Y., Imkeller P. and M. Muller (2005): ``Utility maximization in incomplete markets'', {\it Annals of Applied Probability}, {\bf 15}, 1691-1712. 

\bibitem{jac87} Jacod J. (1987): ``Grosissement initial, hypoth\`ese H' et th\'eor\`eme de Girsanov, S\'eminaire 
de calcul stochastique", Lecture Notes in Maths, {\bf 1118}, 1982-1983. 

\bibitem{jaryu01} Jarrow R.A. and F. Yu (2001): ``Counterparty risk and the pricing of defaultable securities'',
{\it Journal of Finance}, {\bf 56}, 1765-1799.

\bibitem{jealec09} Jeanblanc M. and Y. Le Cam (2009): ``Progressive enlargement of filtrations with initial times'', forthcoming in  {\it Stochastic Processes and their Applications}. 

\bibitem{jeu80} Jeulin T. (1980): ``Semimartingales et grossissements d'une filtration'', Lect. Notes in Maths, Springer, {\bf 883}. 

\bibitem{jeuyor85} Jeulin T. and M. Yor (1985): ``Grossissement de filtration : exemples et applications", Lect. Notes in Maths, Springer, {\bf 1118}. 

\bibitem{jiapha09} Jiao Y. and H. Pham (2009): ``Optimal investment with counterparty risk: a default-density modeling approach'', forthcoming in {\it Finance and Stochastics}.

\bibitem{kaz79} Kazamaki N. (1979): ``A sufficient condition for the uniform integrability of exponential martingales", \textit{Math. Rep. Toyama Univ.}, {\bf 2}, 1Ð11.

\bibitem{kob00} Kobylanski M. (2000): ``Backward stochastic differential equations and partial differential equations with quadratic growth", {\it The Annals of Probability}, {\bf 28}, 558-602. 

\bibitem{limque09} Lim T. and M.C. Quenez (2011): ``Utility maximization in incomplete market with default'', {\it Electronic Journal of Probability}, {\bf16}, 1434-1464.

\bibitem{mor09} Morlais M.A. (2009): ``Utility maximization in a jump market model", {\it Stochastics and Stochastic Reports}, {\bf 81}, 1-27. 

\bibitem{mor10} Morlais M.A. (2010): ``A new existence result for quadratic BSDEs with jumps with
application to the utility maximization problem", {\it Stochastic Processes and  Applications}, 
{\bf 120 (10)}, 1966-1995.

\bibitem{Ok07} Oksendal B. (2007): ``Stochastic Differential Equations, An introduction with Applications", sixth edition, Springer, Berlin. 

\bibitem{pp90} Pardoux E. and S. Peng (1990): ``Adapted solution of a backward stochastic differential equation'', {\it Systems \& Control Letters}, {\bf 37}, 61-74. 



\bibitem{pha09} Pham H. (2010): ``Stochastic control under progressive enlargement of filtrations and applications to multiple defaults risk management",  {\it Stochastic processes and Their Applications},    {\bf 120}, 1795-1820.


\bibitem{pro05} Protter P. (2005): ``Stochastic Integration and Differential Equation'', 2-nd Edition, \textbf{21}, Corrected 3-rd Printing, Stochastic modeling and applied probability, Springer.

 


 

\end{thebibliography}
\end{document}